\documentclass[10pt,twocolumn,twoside,journal]{IEEEtran} 

\usepackage{psfrag}
\usepackage{graphicx}
\usepackage{cite}
\usepackage{amsmath} 
\usepackage{amssymb}  
\usepackage{amsfonts}
\usepackage{mathrsfs}
\usepackage{mathtools}
\usepackage[dvipsnames]{xcolor}
\usepackage[thmmarks, amsmath]{ntheorem}
\usepackage{enumerate}
\usepackage{circuitikz}
\usepackage{booktabs}
\usepackage{colortbl}
\usepackage{gensymb}
\usepackage{pgfplots}
\usetikzlibrary{spy,calc}
\usepackage{etoolbox}
\newbool{phase}
\newbool{mag}
\usepackage{tikz-3dplot}

\DeclareMathOperator{\diag}{diag}

\DeclareMathOperator{\diff}{d}

\providecommand{\norm}[1]{\lVert#1\rVert}
\newcommand{\R}{{\mathbb R}}
\newcommand{\mc}{\mathcal}
\newcommand{\ddt}{\frac{\diff}{\diff t}}

\newtheorem{theorem}{Theorem}

\newtheorem{proposition}{Proposition}
\newtheorem{condition}{Condition}
\newtheorem{assumption}{Assumption}
\newtheorem{remark}{Remark}
\newtheorem{definition}{Definition}

\newcommand{\red}{}

\graphicspath{{Figures/} }

\begin{document}

\title{\LARGE \bf  Global phase and magnitude synchronization of coupled oscillators with application to the control of grid-forming power inverters  \thanks{This work was partially funded by the European Union's Horizon 2020 research and innovation program under grant agreement N$^\circ$ 691800 and the ETH Seed Project SP-ESC 2015-07(4). This article reflects only the authors' views and the European Commission is not responsible for any use that may be made of the information it contains.}}
\author{Marcello Colombino, Dominic Gro\ss, Jean-S\'ebastien Brouillon, and  Florian D\"orfler  \thanks{Marcello Colombino is with the National Renewable Energy Laboratory (NREL). This work was completed before he joined NREL and was not funded by NREL or the US Department of Energy. Dominic Gro\ss, Jean-S\'ebastien Brouillon, and  Florian D\"orfler are with the Automatic Control Laboratory, ETH Z\"urich, Switzerland. e-mail: marcello.colombino@nrel.gov, \{grodo,jeanb,dorfler\}@ethz.ch.}}

\maketitle
\begin{abstract}
In this work, we explore a new approach to synchronization of coupled oscillators. In contrast to the celebrated Kuramoto model we do not work in polar coordinates and do not consider oscillations of fixed magnitude. We propose a synchronizing feedback based on relative state information and local measurements that induces consensus-like dynamics. We show that, under a mild stability condition, the combination of the synchronizing feedback with a decentralized magnitude control law renders the oscillators' almost globally asymptotically stable with respect to set-points for the phase shift, frequency, and magnitude. 
We apply these result to rigorously solve an open problem in control of inverter-based AC power systems. In this context, the proposed control strategy can be implemented using purely local information, induces a grid-forming behavior, and ensures that a network of AC power inverters is almost globally asymptotically stable with respect to a pre-specified solution of the AC power-flow equations. Moreover, we show that the controller exhibits a droop-like behavior around the standard operating point thus making it backward-compatible with the existing power system operation.
\end{abstract}

\section{Introduction}
The electric power grid is undergoing a period of unprecedented change. A major transition is the replacement of bulk generation based on synchronous machines by distributed generation interconnected to the grid via power electronic devices fed by renewable energy sources. This gives rise to scenarios in which either parts of the transmission grid or an islanded distribution grid may operate without conventional bulk generation by synchronous machines. In either case, the power grid faces great challenges due to the loss of the machines rotational inertia and the loss of self-synchronization dynamics inherent to synchronous machines. 

The prevalent approaches to controlling inverters in the future grid are based on mimicking the physical characteristics and controls of synchronous machines~{\cite{zhong2011synchronverters,jouini2016grid,SDA-SJA:13}}. On the one hand, {this approach is appealing because it results in a well-studied closed-loop behavior compatible with the legacy power system. On the other hand, machine-emulation strategies use a system with fast actuation, almost no inherent energy storage, and stringent limits on the output current (the inverter) to mimic a system with slow actuation, significant energy storage, and potentially large peak currents (the generator). Given these vastly different characteristics it is not obvious that machine emulation is a suitable control strategy for grid-forming power inverters \cite{FM-FD-GH-DH-GV:18}. Depending on the implementation emulation-based strategies results in exceedingly complex dynamics that are subject to time-delays in the control and signal processing algorithms \cite{FM-FD-GH-DH-GV:18,GD-TP-PP-XK-FC-XG:15,HB-TI-YM:14}.
}

In a more abstract setting, the problems of synchronization and coordination of coupled systems with limited or no communication have been extensively studied in the control theory literature~\cite{sepulchre2007stabilization,jadbabaie2003coordination,cucker2007emergent,stan2007analysis} as they find application in many domains beyond power systems, ranging from physics and biology~\cite{strogatz2000kuramoto,strogatz1993coupled,hioe1978quantum,ford1965statistical} to spacecraft coordination~\cite{beard2001coordination,smith2007closed}, and they pose interesting theoretical challenges.

The connection between coupled oscillators and control of power inverters has been extensively studied in the literature. Most of the standard approaches consider so-called droop control \cite{MCC-DMD-RA:93,JMG-MC-TLL-PLC:13a} and, for the purpose of analysis, rely on modified versions of the Kuramoto model of coupled oscillators~\cite{simpson2013synchronization,dorfler2013synchronization,dorfler2012synchronization}. While providing useful insights, these approaches often neglect the voltage dynamics, the associated phasor models are well-defined only near the synchronous steady-state, and the synchronization guarantees are only local as the phasor dynamics admit multiple equilibria on the torus, corresponding to different relative angle configurations. More recent approaches rely on controlling the inverters to behave like virtual Li\'enard-type oscillators~{\cite{johnson2014synchronization,johnson2016synthesizing,MS-FD-BJ-SD:14b,sinha2016synchronization,LABT-JPH-JM:12,LABT-JPH-JM:13}}. This approach is promising, because virtual oscillator control can globally synchronize an inverter-based power system and has been validated experimentally. 
{However, in its original form virtual oscillator control (VOC) cannot be dispatched to track references for the power injection and voltage. Likewise, all theoretic investigations are limited to synchronization with identical angles and voltage magnitudes, i.e., no power is exchanged between the inverters. For passive loads it can be shown that power is delivered to the loads \cite{LABT-JPH-JM:13}, but the power sharing by the inverters and their voltage magnitudes are determined by the load and network parameters. As of today, it is unclear how to extend the basic virtual oscillator controller to track references for the power injection, i.e., to steer the system to a desired solution of the power flow equations with nonzero relative voltage angles and magnitudes between the inverters. However, to operate grid-forming inverters as part of a future grid they need to be dispatchable, i.e., able to track references for the power injection and voltage.}

To overcome these challenges of emulation, droop, and virtual oscillator control, we follow a top-down approach inspired by consensus strategies and flocking~ \cite{fax2004information,ren2008distributed,lin2013algorithms,sepulchre2007stabilization,jadbabaie2003coordination,cucker2007emergent}. {We analyze a novel control strategy that was recently proposed by the authors in the context of power systems \cite{colombino2017global}. In particular, it is established in \cite{colombino2017global} that there exists set-points for the relative phase angles that ensure almost global synchronization of a network of power inverters. The key result of the present manuscript is an insightful stability condition that quantifies stability in the parameter space, i.e., we provide a tractable condition on the set-points, control gains, and network parameters that ensures almost global asymptotic stability of the proposed controller. Moreover, we provide a rigorous and deep theoretical analysis and link and compare the results to the literature on synchronization of coupled nonlinear oscillators. We show that, under a mild stability condition, the proposed controller renders the oscillators' almost globally asymptotically stable with respect to pre-specified set-points for the phase shift, frequency, and magnitude. In the context of networked inverters, we show how the underlying network physics can be exploited to implement the controller in a decentralized fashion, i.e., only using local measurements and set-points for power injection and voltage magnitude. To the best of our knowledge, the proposed solution is the first control strategy that, under standard assumptions, ensures almost global asymptotic stability of an inverter-based AC power system with respect to a pre-specified solution of the AC power flow equations. Finally, we show that the proposed controller, much like virtual oscillator control and most machine-emulation strategies, behaves similarly to a conventional droop controller around the synchronous harmonic steady state.
}
%
%
%
%

The remainder of this section recalls some basic notation and results from graph theory. In Section~\ref{sec.control.problem} we introduce the class of oscillators under consideration, we formally state the control objectives, and we propose a synchronizing closed-loop behavior that admits a decentralized implementation. Section~\ref{sec.inverter.model} introduces the model of an inverter based power system in the context of the problem setup presented in Section~\ref{sec.control.problem} as well as a reformulation of the proposed controller in the context of power systems.  The main result, in which we prove almost global {asymptotic stability}, is given in Section~\ref{sec.proofs}. Section~\ref{sec.example} presents a simulation example, and Section~\ref{sec.conclusion} provides the conclusion and outlook.

\subsection*{Notation}
The set of real numbers is denoted by $\R$. ${\R_{>0}}$ denotes the set of nonnegative reals{, and $\mathbb{S}^1 = [0,2\pi]$ denotes the circle.} Given $\theta\in{\mathbb{S}^1}$ we define
\begin{align*}
R(\theta) \coloneqq \begin{bmatrix}\cos(\theta) & -\sin(\theta)\\  \sin(\theta) & \cos(\theta) \end{bmatrix},~ \text{and}~ J\coloneqq \begin{bmatrix} 0 & -1\\  1 & 0 \end{bmatrix}.
\end{align*}
Given a matrix $A$, $A^\top$ denotes its transpose. We use $\bar{\sigma}(A)$ to indicate the largest singular value of $A$. We write $A\succcurlyeq0$  $(A\succ0)$ to denote that $A$ is symmetric and positive semidefinite (definite). {For column vectors $x\in\R^n$ and $y\in\R^m$ we use $(x,y) = [x^\top, y^\top]^\top \in \R^{n+m}$ to denote a stacked vector.} Furthermore, $I_n$ denotes the identity matrix of dimension $n$, $\otimes$ denotes the Kronecker product, $\norm{x}$ denotes the Euclidean norm, and we use $\norm{x}_C \coloneqq \min_{z \in \mathcal{C}} \norm{z-x}$ to denote the distance of a point $x$ to a set $\mathcal{C}$. We use $\varphi_f(t,x_0)$ to denote the solution of $\ddt x = f(x)$ at time $t \geq 0$ starting from the initial condition $x(0)=x_0$ at time $t_0=0$. Moreover, we use ${\mathscr{B}_{\mc C}(r)} \coloneqq \{ x \in \mathbb{R}^n \vert \norm{v}_{\mc C} {<} r \}$ to {denote} an $r$-neighborhood of a set $\mc C$. Finally, given $x\in\R$, we denote by $\lceil x \rceil\coloneqq \min_{y\in\mathbb N, y\ge x} y$ the ceiling operator that rounds up a real number to the nearest integer.

\subsection*{Preliminaries from graph theory}
In this paper, we will make extensive use of concept from algebraic graph theory. We refer the reader to~\cite{lns-v.95} for a comprehensive review of the topic. {Given a simple connected graph (i.e., an undirected graph containing no self-loops or multiple edges) denoted by} $\mc G = (\mc V, \mc E,{\{w_{jk}\}_{(j,k)\in \mc E}})$, $\mc V = \{1,..., N\}$ is the set of nodes{,  $\mc E \subset {(\mc V \times \mc V) \setminus \cup_{i \in \mc V} (i,i)}$, with $|\mc E|=M$, is the set of undirected edges, and $\{w_{jk}\}_{(j,k)\in \mc E}$ is the set of edge weights $w_{jk}$ corresponding to the respective edge. Assigning an index $\ell \in \{1, . . . ,M\}$ and an arbitrary direction to each edge $(i, j) \in \mc E$, the (oriented) incidence matrix $B \in \R^{N\times M}$ is defined component-wise by $B_{k\ell} = 1$ if node $k$ is the sink node of edge $\ell$ and by $B_{k\ell} = -1$ if node $k$ is the source node of edge $\ell$; all other elements are zero. Given $x \in \R ^N$, the vector $B^\top x$ has components $x_i - x_j$ corresponding to the oriented edge from $j$ to $i$. Using the diagonal matrix of edge weights $\diag\left (\{w_\ell\}_{\ell\in\mc E} \right)$, the graph Laplacian  $L$ is defined as 
\begin{align}\label{eq.laplacian.incidence}
L \coloneqq B\,  \diag(\{w_\ell\}_{\ell \in \mc E} )B^\top.
\end{align}
}

\section{Problem setup and control strategy}\label{sec.control.problem}
Before focusing on the main engineering application, we will first precisely state the problem at hand as an abstract coordination problem that arises in many scientific disciplines \cite{sepulchre2007stabilization,jadbabaie2003coordination,cucker2007emergent}. Specifically, consider $N$ two-dimensional integrators with a state variable $v_k\in\R^2$ {that} represents the state of oscillator {$k \in \{1,\ldots,N\}$}. {In a polar coordinate representation, $\norm{v_k}$ is the magnitude of the oscillation, and the angle $\theta_k \in \mathbb{S}^1$ associated to each oscillator satisfies $v_k = R(\theta_k) (\norm{v_k},0)$.} The oscillators are interconnected via a simple connected {(not necessarily complete)} graph $\mc G = (\mc V, \mc E,{\{w_{jk}\}_{(j,k)\in \mc E}})$, with vertices {$\mc V$ corresponding to the oscillators and} Laplacian matrix $L$. {It should also be noted that the graph $\mc G$ is not necessarily modeling a communication network. For example, in the context of power systems the graph models the power transmission network.} We define the extended Laplacian as $\mc L \coloneqq L\otimes I_2$. To each vertex we associate a two-dimensional integrator, a local controller $u_k(v_k,y_k): \R^2 \times \R^2 \to \R^2$, and an output variable $y_k\in\R^2$. 
The dynamics of the closed loop are given by
\begin{subequations}\label{eq.oscillators}
\begin{align}
\frac{\diff}{\diff t} v_k &  = u_k(v_k,y_k)\quad k\in\mc V,\\
y & =  \mc L v,
\end{align}
\end{subequations}
where $v=(v_1,\ldots,v_N) \in\R^{2N}$ is the overall state vector, and $y=(y_1,\ldots,y_N) \in\R^{2N}$ is a vector of relative outputs. {We aim to construct decentralized control laws $u_k(v_k,y_k)$ that synchronize the $N$ two-dimensional integrators (almost) globally to a harmonic oscillations with prescribed frequency, relative phase angles, and magnitudes. Following the definition commonly used in control of large-scale systems \cite{lunze1992feedback} we consider a controller to be decentralized if it only uses the local state $v_k$, the local output $y_k$, and set-points that may be updated infrequently by a centralized decision maker. This is also the setup encountered in the power system application discussed in Section \ref{sec.inverter.model}\red{, in which $v_k$ corresponds to the AC voltage of the $k$-th inverter.}}
\subsection{Synchronization: a geometric characterization} \label{sec.control.problem.obj}
We now give a geometric characterization of the control objectives that allows us to formalize the objective of synchronization {and magnitude regulation}, and re-state it as the problem of stabilizing the oscillators with respect to a set. Let us define ${\mc J\coloneqq  I_N\otimes J}$. Given a function $f:\R^{2N\times 2N} \to \R^{2N\times 2N}$ that describes the closed-loop dynamics $\frac{\diff}{\diff t} v = f(v)$, a frequency $\omega_0\ge0$, relative phase set-points $\{{\theta_{k1}^\star}\}_{k=2}^N$ with $\theta_{k1}\in {\mathbb{S}^1}$, and magnitude set-points $v^\star$ we define sets corresponding to the three steady-state specifications:
\begin{itemize}
\item \textbf{Relative phase {angles}:}
\begin{align}\label{eq.obj.theta}
\mc S \coloneqq \left\{ v \in \R^{2N} \left\vert\, \frac{v_k}{{v^\star_k}} -  R(\theta_{k1}^\star)\frac{v_1}{v^\star_1} = 0, \, \forall k \in \mc V \setminus \{1\} \right.\right\},
\end{align}
\item \textbf{{Synchronous frequency:}}
\begin{align}\label{eq.obj.sync}
\begin{split}
&\mc F_{f,\omega_0} \coloneqq  \left\{ v \in \R^{2N}  \left\vert\; f(v) = \omega_0 \mc J v\right.\right\},
\end{split}
\end{align}
\item \textbf{Magnitude {of the oscillations}:}
\begin{align}\label{eq.obj.norm}
\mc A_k \coloneqq \left\{ v_k \in \R^2 \left\vert\; \norm{v_k} = {v^\star_k} \right.\right\}.
\end{align}
\end{itemize}
Finally, we define the product set $\mc A \coloneqq \times^{N}_{k=1} \mc A_k$. In $\mc S$ the oscillators will be at the prescribed phase shift with respect to each other, in $\mc F_{f,\omega_0}$ they rotate at frequency $\omega_0$ and behave like harmonic oscillators and in $\mc A$ their state variables have the correct magnitudes. The aim of this work is to construct a vector field $f$ of the form
\[
f(v)=\begin{bmatrix}
u_1(v_1, (\mc L v)_1)\\
\vdots\\
u_N(v_N, (\mc L v)_N)
\end{bmatrix},
\]
such that the intersection of the three sets, i.e., 
\begin{align}
 \mc T_{f,\omega_0} \coloneqq \mc S \cap  \mc F_{f,\omega_0} \cap \mc A,
\end{align}
is invariant and (almost) globally asymptotically stable.

\subsection{Synchronizing controllers and closed-loop dynamics}\label{sec.consensus.control}
The challenge of controlling coupled oscillators is that each controller $u_k(v_k,y_k)$ can rely only on measurements of the local state and output variables,  while the control objectives are global. In this section, we introduce a decentralized globally synchronizing control law for the oscillators. Before doing that we need some preliminary definitions. Given the graph $\mc G$ and angle set-points $\{{\theta_{k1}^\star}\}_{k=1}^N$ with $\theta^\star_{k1}\in {\mathbb{S}^1}$ we define $\theta_{jk}^\star\coloneqq \theta_{j1}^\star - \theta_{k1}^\star$, the matrices 
\begin{align}\label{eq.Kk.def}
K_k  \coloneqq  \sum\nolimits_{j:(j,k)\in\mc E}w_{jk}\left(I_2-\frac{v^\star_j}{v^\star_k}R(\theta_{jk}^\star)\right), \quad k\in\mc V,
\end{align}
as well as the functions 
\begin{align}\label{eq.phi}
\Phi_k(v_k)  \coloneqq   \frac{{v^\star_k}-\norm{v_k}}{{v^\star_k}},
\end{align}
{which compute a scaled magnitude error.}
We propose the following control laws $u_k(v_k,y_k)$ 
\begin{align} \label{eq.controller}
 \frac{\diff}{\diff t} v_k  & =   u_k(v_k,y_k) 
 \nonumber\\
 & = \omega_0 Jv_k   +  {\eta  \left(  K_k v_k - y_k \right) + \alpha \Phi_k(v_k)I_2 v_k },
\end{align}
 where $\eta>0$ {is the synchronization gain and $\alpha>0$ is the control gain associated with amplitude regulation}. For brevity of the presentation, we consider the case in which all oscillators have the same control gains. The decentralized controllers~\eqref{eq.controller}, together with the fact that $y = \mc L v$ give rise to the following closed loop
\begin{align}\label{eq.desired.cl.magnitude.two}
 \frac{\diff}{\diff t} v =  \left(\omega_0\mc J  + {\eta} (\mc K - \mc L)  + \alpha \Phi(v) \right)v {\; \eqqcolon  f(v)},
\end{align}
where $\mc K$ and $\Phi(v)$ are defined as $\mc K \coloneqq \diag(\{K_k\}_{k=1}^N)$ and $\Phi(v)  \coloneqq   \diag\left(\{\Phi_k(v_k)I_2\}_{k=1}^N\right)$.

We can give an intuitive interpretation for the closed-loop dynamics~\eqref{eq.desired.cl.magnitude.two} by defining the phase error $e_\theta(v) \coloneqq [ e^\top_{\theta,1}(v), \ldots, e^\top_{\theta,N}(v)]^\top$ between each oscillator and its neighbors
\begin{align}\label{eq.phase.err}
  e_{\theta,k}(v) \coloneqq \sum\nolimits_{j:(j,k)\in\mc E} w_{jk} \left(v_j - \frac{v^\star_j}{v^\star_k}R(\theta_{jk}^\star) v_k \right ),
\end{align}
and the {magnitude error} as $e_{\norm{v}}(v)\coloneqq[ e^\top_{\norm{v},1}(v), \ldots, e^\top_{\norm{v},N}(v)]^\top$, where the components $e_{\norm{v},k}(v)$ are defined as
\begin{align}\label{eq.error.phi}
 e_{\norm{v},k}(v)\coloneqq\Phi_k(v_k)v_k.
\end{align}
While the magnitude error depends only {on the local oscillator state}, the phase error depends on the neighbors' state variables. {To clarify the interpretation of the control law~\eqref{eq.controller}, the following proposition restates the closed-loop dynamics~\eqref{eq.desired.cl.magnitude.two} in terms of phase and magnitude errors.}
\begin{proposition}{{\bf(Interpretation of the controller)}}\\
The dynamics~\eqref{eq.desired.cl.magnitude.two} are equivalent to
\begin{align}\label{eq.desired.cl.magnitude}
 \frac{\diff}{\diff t} v =  \omega_0\mc J v + \eta e_\theta(v) + \alpha e_{\norm{v}}(v) .
\end{align}
\end{proposition}
\begin{IEEEproof}
We note that:
\begin{align}\label{eq.error.theta}
\begin{split}
 &e_{\theta,k}(v) = \!\sum_{j:(j,k)\in\mc E} \!\!w_{jk} \left( v_j-\frac{v^\star_j}{v^\star_k}R(\theta_{jk}^\star) v_k - v_k + v_k \right) \\
 &=\!\sum_{j:(j,k)\in\mc E} \!\!w_{jk} (v_j\!-v_k) + \!\sum_{j:(j,k)\in\mc E}\!\!w_{jk}\left(I_2\!-\frac{v^\star_j}{v^\star_k}R(\theta_{jk}^\star)\right) v_k, \\
 & = - y_k +  K_k  v_k.
 \end{split}
\end{align}
Form ~\eqref{eq.error.theta} and the fact that $y=\mc L v$ we conclude that
\begin{align}\label{eq.error.theta2}
\begin{split}
 &e_{\theta}(v) =  \mc K  v - y = (\mc K-\mc L)v.
 \end{split}
\end{align}
In view of~\eqref{eq.error.theta2} and~\eqref{eq.error.phi} the closed loop~\eqref{eq.desired.cl.magnitude.two} is equivalent to~\eqref{eq.desired.cl.magnitude} and the proof is complete.
\end{IEEEproof}
{In words, the control law \eqref{eq.controller} induces closed-loop dynamics that are a combination of a harmonic oscillation with frequency $\omega_0$ as well as a linear feedback of the phase error $e_\theta(v)$ and magnitude error $e_{\norm{v}}(v)$.
The corresponding components of $u_k$ are illustrated in Figure~\ref{fig.vector.field}: the first term in~\eqref{eq.desired.cl.magnitude} is tangential to $v_k$ and corresponds to the rotational speed / frequency $\omega_0$, the second term $\eta e_\theta$ synchronizes the relative phase angles $\theta_{jk}$ to the set-points $\theta_{jk}^\star$, and the third term $\alpha e_{\norm{v}}$ results in a component that is radial to $v_k$ and regulates the amplitude.} Note that, if the phase angles and magnitudes are correct, i.e., $ e_\theta= { e_{\norm{v}}}  = 0$, the residual dynamics of~\eqref{eq.desired.cl.magnitude} are those of $N$ decoupled harmonic oscillators with synchronous frequency $\omega_0$. This gives an intuition on the synchronizing behavior of the controller~\eqref{eq.controller}. In Theorem~\ref{thm.ags} (Section~\ref{sec.proofs}) we prove that, under a mild stability condition,~\eqref{eq.controller} drives the system to synchrony and satisfies all objectives described in Section~\ref{sec.control.problem.obj} from almost all initial conditions. 
\booltrue{phase}
\booltrue{mag}
\pgfmathsetmacro{\rvec}{0.5}
\pgfmathsetmacro{\valeta}{0.5}
\pgfmathsetmacro{\valapha}{2}
\pgfmathsetmacro{\valomega}{0.5}
\pgfmathsetmacro{\valtheta}{30}
\pgfmathsetmacro{\etheta}{1}
\pgfmathsetmacro{\emag}{1}
\pgfmathsetmacro{\vonex}{0.7*\rvec*cos(-5)}
\pgfmathsetmacro{\voney}{0.7*\rvec*sin(-5)}
\pgfmathsetmacro{\vtwox}{1.2*\rvec*cos(65)}
\pgfmathsetmacro{\vtwoy}{1.2*\rvec*sin(65)}
\begin{figure}[htbp]
\begin{center}
\tdplotsetmaincoords{0}{0}
\pgfmathsetmacro{\sintheta}{sin(\valtheta)}%
\pgfmathsetmacro{\costheta}{cos(\valtheta)}%

\pgfmathsetmacro{\normvone}{sqrt(\vonex*\vonex+\voney*\voney)}%
\pgfmathsetmacro{\normvtwo}{sqrt(\vtwox*\vtwox+\vtwoy*\vtwoy)}%

\begin{tikzpicture}[scale=4,tdplot_main_coords]

\node (O) at (0,0,0) [anchor=north]  {$0$};



\draw[thick,ForestGreen!50] (\rvec,0,0) arc  (0:360:\rvec);

\draw[-o, very thick,NavyBlue, shorten >=-1.15mm] (0,0) -- (\vonex,\voney) node[anchor=north west]{ $v_1$};
\draw[-o, very thick,NavyBlue, shorten >=-1.15mm] (0,0) -- (\vtwox,\vtwoy) node[anchor=north west]{ $v_2$};

\draw[->, very thick,ProcessBlue]  (\vonex,\voney) -- (\vonex-\valomega*\voney,\voney+\valomega*\vonex) node[anchor= south]{ $\omega_0 J\, v_1$};
\draw[->, very thick,ProcessBlue]  (\vtwox,\vtwoy)  -- (\vtwox-\valomega*\vtwoy,\vtwoy+\valomega*\vtwox) node[anchor=east ]{ $\omega_0 J\, v_2$};

\ifbool{phase}{
\draw[->, very thick,orange]  (\vonex,\voney) -- (\vonex+\valeta*\vtwox-\valeta*\costheta*\vonex+\valeta*\sintheta*\voney,
\voney+\valeta*\vtwoy-\valeta*\sintheta*\vonex-\valeta*\costheta*\voney) node[anchor=east ]{ $e_{\theta^\star,1}$};

\draw[->, very thick,orange]  (\vtwox,\vtwoy) -- (\vtwox+\valeta*\vonex-\valeta*\costheta*\vtwox-\valeta*\sintheta*\vtwoy,
\vtwoy+\valeta*\voney+\valeta*\sintheta*\vtwox-\valeta*\costheta*\vtwoy) node[anchor=west ]{ $e_{\theta^\star,2}$};
}

\ifbool{mag}{
\draw[->, very thick,Orchid]  (\vonex,\voney) -- (\vonex+\valapha*\rvec*\vonex-\valapha*\vonex*\normvone,
\voney+\valapha*\rvec*\voney-\valapha*\voney*\normvone) node[anchor=south]{ $e_{v^\star,1}$};

\draw[->, very thick,Orchid]  (\vtwox,\vtwoy) -- (\vtwox+\valapha*\rvec*\vtwox-\valapha*\vtwox*\normvtwo,
\vtwoy+\valapha*\rvec*\vtwoy-\valapha*\vtwoy*\normvtwo) node[anchor= east]{ $e_{v^\star,2}$};
}

\end{tikzpicture}
\caption{Illustration of the components of the closed-loop vector field for two interconnected oscillators.}
\label{fig.vector.field}
\end{center}
\end{figure}
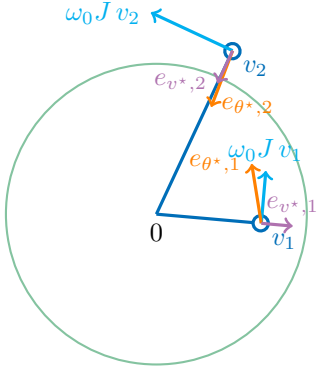

\subsection{Relationship to the Kuramoto oscillator}\label{sec.kuramoto}
The control objectives outlined in Section~\ref{sec.control.problem.obj} require that each oscillator converges to the circle defined by the set $\mc A_k$ and synchronizes with the other oscillators to a rotating trajectory with a given frequency and prescribed angle differences to neighboring oscillators. This problem falls in the broader class of coordination/synchronization problems using neighbors-only communication extensively studied in the literature ~\cite{sepulchre2007stabilization,jadbabaie2003coordination,cucker2007emergent}. In the field of coupled oscillators, most of the literature is focused on solving such synchronization problems while restricting the dynamics of each system to evolve on the circle $\mc A_k$ at all times. The most celebrated example of such oscillator is the Kuramoto model~\cite{kuramoto1975self,acebron2005kuramoto,REM-SHS:05,FD-FB:13b}:
\begin{align}\label{eq.kuramoto}
\begin{split}
\dot \theta_k = \omega_k +  {\eta}\sum\nolimits_{j:(j,k)\in\mc E} w_{jk} \sin(\theta_j-\theta_k),\quad k\in\mc V
 \end{split}
\end{align}
The linear part of the closed-loop dynamics in~\eqref{eq.desired.cl.magnitude.two} has been described in the literature as the ``linear reformulation of the Kuramoto model"~\cite{roberts2008linear}.  With the following proposition, we show that the Kuramoto model~\eqref{eq.kuramoto} can be derived by restricting the evolution of the dynamics  \eqref{eq.desired.cl.magnitude.two} to the unit circle (i.e., $\mc A_k$ with ${v^\star_k} = 1$) by a projection operation. We define the projector onto {the direction orthogonal to}  $v_k$ as
\[
\Pi_{{v^\perp_k}} = I_2-\frac{v_kv_k^\top}{\norm{v_k}^2},
\]
Moreover, we define $\Pi_{{v^\perp}} \coloneqq  \diag\left(\{\Pi_{{v^\perp_k}}\}_{k=1}^N\right)$
\begin{proposition}{\bf{(Relationship to the Kuramoto model)}}\\
Consider the vector field $f(v)$ defined in~\eqref{eq.desired.cl.magnitude.two}, {the angle $\theta_k \in \mathbb{S}^1$ such that $v_k = R(\theta_k) ({\norm{v_k}},0)$}, and ${v^\star_k}=1$. For all $ v(0) \in \mc A$ the dynamical system
\begin{align}\label{eq.kuramoto.equivalent}
 \frac{\diff}{\diff t} v = \Pi_{{v^\perp}} f(v)
\end{align}
is equivalent to the Kuramoto model of coupled oscillators~\eqref{eq.kuramoto} {with natural frequencies given by $\omega_k\coloneqq \omega_0 + \eta\sum\nolimits_{j:(j,k)\in\mc E}w_{jk} \sin(\theta^\star_{jk}), \quad k\in\mc V$}.

\end{proposition}
\begin{IEEEproof}
{Note that, {by definition of} the projection operator $\Pi_{{v^\perp}}$, $\cal A$ is invariant under the dynamics~\eqref{eq.kuramoto.equivalent}, i.e., $\norm{v_k(t)}=1$ for all $t>0$.} Consider the $k^\text{th}$ (2-dimensional) component of~\eqref{eq.kuramoto.equivalent}:
\begin{align}\label{eq.Kuramoto.equivalent.steps}
\begin{split}
 \frac{\diff}{\diff t} v_k & =  \omega_0 Jv_k  + \\
 &+ \eta \,{\Pi_{{v^\perp_k}} \left(  K_k v_k + \sum\nolimits_{j:(j,k)\in\mc E} w_{jk} (v_j\!-v_k) \right)}\\
 &+\alpha \,{\Pi_{{v^\perp_k}}\left( \Phi_k(v_k)I_2  v_k\right)} \\ 
  & = \omega_k Jv_k +  {\eta} {\Pi_{{v^\perp_k}}\left(  \sum\nolimits_{j:(j,k)\in\mc E} w_{jk} v_j\right)}. \\
 \end{split}
\end{align}
{To prove the last equality in~\eqref{eq.Kuramoto.equivalent.steps}, note that with $v^\star_k=1$ for all $k\in\mc V$, $K_k = \sum\nolimits_{j:(j,k)\in\mc E} w_{jk} \cos(\theta_j-\theta_k)I_2 + \sum\nolimits_{j:(j,k)\in\mc E} w_{jk} \sin(\theta_j-\theta_k)J$. This, together with the fact that $\|v_k\|=1$ implies that $\Pi_{{v^\perp_k}} K_k v_k = \sum\nolimits_{j:(j,k)\in\mc E} w_{jk} \sin(\theta_j-\theta_k)Jv_k$}.
The dynamics~\eqref{eq.Kuramoto.equivalent.steps}, written in polar coordinates following the specifications outlined in~\cite[Proposition 2]{zhu2013synchronization}, reduce to~\eqref{eq.kuramoto}.
\end{IEEEproof}
The Kuramoto oscillator has been extensively studied, and the literature provides numerous conditions that guarantee local synchronization~\cite{FD-FB:13b}. {Observe that a large amplitude control gain $\alpha$ implies that $\norm{v_k} \approx v^\star_k$, i.e., \eqref{eq.kuramoto.equivalent} is the limiting vector field of \eqref{eq.desired.cl.magnitude.two} for high gain amplitude control. However,} while closely related to the Kuramoto oscillator, the closed loop~\eqref{eq.desired.cl.magnitude.two} does not fix the magnitude of the state variables during transients and, {as outlined in Section~\ref{sec.proofs}}, converges to the torus asymptotically {for a small amplitude control gain $\alpha$}. This will allow us to prove conditions that guarantee almost global instead of local synchronization.

\section{Defining engineering application: inverter-based power grids}\label{sec.inverter.model}
In the following, we motivate the abstract theoretical problem setting in Section~\ref{sec.control.problem} by connecting it to the control of inverter based AC power systems. The main challenge of controlling power inverters in a power system that is possibly islanded and without synchronous machines is that each controller can only rely on local measurements, while the control objectives (such as frequency synchronization and convergence to a desired solution of the power-flow equations) are global in nature. In this section, we show how the control of an inverter-based power grid such as the one illustrated in Figure~\ref{fig.inverters}, can be abstractly modeled within the setup described in Section~\ref{sec.control.problem.obj}.

\subsection{Modelling of inverter-based power grids}
\begin{figure}[htbp]
\begin{center}
\begin{circuitikz}[american voltages]

\ctikzset{bipoles/resistor/height=0.15}
\ctikzset{bipoles/resistor/width=0.4}

\ctikzset{bipoles/generic/height=0.15}
\ctikzset{bipoles/generic/width=0.4}

\ctikzset{bipoles/length=.6cm}

\coordinate (I1) at (0,0);

\draw ($ (I1) + (0,-1) $) node [ground] (g1) {};
\draw (g1) to[sV, l_=\mbox{$\frac{\diff}{\diff t} {{{v}}_{{k}}(t)}= u_k(v_k,R(\kappa)i_{o,k})$}] (I1);
\draw ($ (I1) + (-0.5,-1) $) node [ground] (g12) {};
\draw (g12) to[R] ($ (I1) + (-0.5,-0.3) $);
\draw (I1) to[short] ($ (I1) + (-0.5,-0.3) $);
\draw (I1) to[short,i^=\mbox{$i_{o,k}$}] ($ (I1) + (2,0) $);
\draw [dashed] ($ (I1) + (2,0) $) to[short,l={to network}] ($ (I1) + (3,0) $);
\end{circuitikz}
\caption{Schematics of the decentralized control setup: the voltage $v_k$ is fully controllable but only local measurements can be used.}
\label{fig.inverter.controlled}
\end{center}
\end{figure}
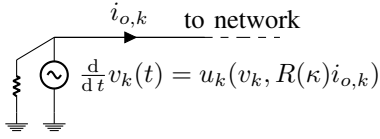

\begin{figure}[htbp]
\begin{center}
\begin{circuitikz}[american voltages]

\ctikzset{bipoles/resistor/height=0.15}
\ctikzset{bipoles/resistor/width=0.4}

\ctikzset{bipoles/generic/height=0.15}
\ctikzset{bipoles/generic/width=0.4}

\ctikzset{bipoles/length=.6cm}

\coordinate (I1) at (0,0);
\coordinate (I2) at (4,0.5);
\coordinate (I3) at (2,-1);
\coordinate (I4) at (1,1.5);

\draw ($ (I1) + (0,-1) $) node [ground] (g1) {};
\draw (g1) to[sV] (I1);
\draw ($ (I1) + (-0.5,-1) $) node [ground] (g12) {};
\draw (g12) to[R] ($ (I1) + (-0.5,-0.3) $);
\draw (I1) to[short] ($ (I1) + (-0.5,-0.3) $);

\draw ($ (I2) + (0,-1) $) node [ground] (g2) {};
\draw (g2) to[sV] (I2);
\draw ($ (I2) + (0.5,-1) $) node [ground] (g22) {};
\draw (g22) to[R] ($ (I2) + (0.5,-0.3) $);
\draw (I2) to[short] ($ (I2) + (0.5,-0.3) $);

\draw ($ (I3) + (0,-1) $) node [ground] (g3) {};
\draw (g3) to[sV] (I3);
\draw ($ (I3) + (0.5,-1) $) node [ground] (g32) {};
\draw (g32) to[R] ($ (I3) + (0.5,-0.3) $);
\draw (I3) to[short] ($ (I3) + (0.5,-0.3) $);

\draw ($ (I4) + (0,-1) $) node [ground] (g4) {};
\draw (g4) to[sV] (I4);
\draw ($ (I4) + (0.5,-1) $) node [ground] (g42) {};
\draw (g42) to[R] ($ (I4) + (0.5,-0.3) $);
\draw (I4) to[short] ($ (I4) + (0.5,-0.3) $);

\node (I12) at ($(I1)!0.5!(I2)$) {};
\draw 		(I1) 
to[R,*-] (I12)
to [L,-*] (I2);

\node (I13) at ($(I1)!0.5!(I3)$) {};
\draw 		(I1) 
to[R,*-] (I13)
to [L,-*] (I3);

\node (I32) at ($(I3)!0.5!(I2)$) {};
\draw 		(I3) 
to[R,*-] (I32)
to [L,-*] (I2);

\node (I14) at ($(I1)!0.5!(I4)$) {};
\draw 		(I1) 
to[R,*-] (I14)
to [L,-*] (I4);

\node (I24) at ($(I2)!0.5!(I4)$) {};
\draw 		(I2) 
to[R,*-] (I24)
to [L,-*] (I4);

\end{circuitikz}
\caption{Inverter-based grid}
\label{fig.inverters}
\end{center}
\end{figure}
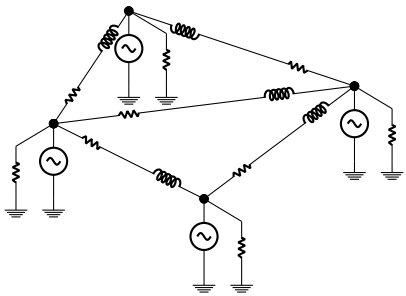

We study the control of $N$ three-phase inverters interconnected by a resistive-inductive network as in Figure~\ref{fig.inverters} with a setup similar to~\cite{ johnson2014synchronization,johnson2016synthesizing,MS-FD-BJ-SD:14b,sinha2016synchronization,LABT-JPH-JM:13}. Each inverter is abstracted as a {three-phase controllable voltage source. Moreover, to each line we associate a three-phase current. All electrical quantities in the three-phase network are assumed to be balanced. Therefore, by applying the well-known Clarke transformation to the three-phase variables we can work in $\alpha\beta$ coordinate frame using two-dimensional vectors for the voltages and currents~\cite{clarke1943circuit}. Appendix~\ref{app.clarke} reviews of the Clarke transformation. To each inverter, we associate a controllable voltage $v_k = (v_{\alpha,k},v_{\beta,k}) \in\R^2$ and an output current $i_{o,k}\in\R^2$. The quantities are illustrated in Figure~\ref{fig.inverter.controlled}.} 

We model the inverter-based power grid as a {simple}, connected, and undirected graph $\mc G~=~(\mc V, \mc E,{\{w_{jk}\}_{(j,k)\in \mc E}})$, where $\mc V = \{1,..., N\}$ is the set of nodes corresponding to the inverters,  $\mc E \subset \mc V \times \mc V$, with $|\mc E|=M$, is the set of undirected edges corresponding to the transmission lines, and ${\{w_{jk}\}_{(j,k)\in \mc E}}$ is the {set of weights} $w_{jk}$ defined as 
\begin{align}\label{eq.weights}
w_{jk}  \coloneqq \frac{1}{\sqrt{r_{jk}^2+\omega_0^2\ell_{jk}^2}},
\end{align}
where $r_{jk}$ and $\ell_{jk}$ are respectively the resistance and inductance of the line $(j,k) \in \mc E$. 

Given the {oriented} incidence matrix {$B$} of the graph ${\mc G}$, we can define $\mc B \coloneqq B\otimes I_2$ that, again, duplicates each edge for the $\alpha$ and $\beta$ components.  We construct the diagonal matrices $R_T\in\R^{M\times M}$ and $L_T\in\R^{M\times M}$ that contain the values of the resistance and inductance of every line on the main diagonal and the matrices $\mathbf R_T \coloneqq R_T\otimes I_2$ and $\mathbf L_T \coloneqq L_T\otimes I_2$ that duplicate the inductance and resistance values for the $\alpha$ and $\beta$ components. 

We assume that the inverters are controlled such that the time constants of the transmission lines are much faster than those of \red{the terminal voltages.} Therefore, we assume that the lines reach steady state almost instantaneously and couple all inverters through the algebraic relationship
\begin{align}\label{eq:quasisteady}
i_o = \mc Y v,
\end{align}
where $i_o = [i_{o,1}^\top,...,i_{o,N}^\top]^\top\!$, $v = [v_1^\top,...,v_N^\top]^\top\!$, and $\mc Y$ is the {admittance} matrix of the network, which can be constructed as
\begin{align}\label{eq.impedance}
\mc Y = \mc B \left(\mathbf R_T + \omega_0\mathcal J \mathbf L_T\right)^{-1} \mc B^\top.
\end{align}
Note that, since $\alpha\beta$ coordinates can be thought as an embedding of the complex numbers into real-valued Euclidean coordinates, the $90^\circ$ rotation matrix $J$ plays the same role that the imaginary unit $\sqrt{-1}$ plays in complex coordinates.
\red{
\begin{remark}
While the quasi-steady-state approximation \eqref{eq:quasisteady} (also known as phasor approximation) of the transmission network is typically justified for conventional power systems, the electromagnetic transients of the transmission lines can compromise the stability of an inverter-based power system \cite{MV+05,VHA+17}. Using results from singular perturbation analysis~\cite{K02} our analysis can be extended to obtain stability guarantees that explicitly include the dynamics of transmission lines. These results are outside of the scope of this work and can be found in \cite{GCB+18}.
\end{remark}
Next,} note the similarity between the expression for a graph's Laplacian matrix~\eqref{eq.laplacian.incidence} and the impedance matrix of the power network~\eqref{eq.impedance}. To formalize this observation, let us make the following standard assumption.
\begin{assumption}{\bf(Uniform inductance-resistance ratio)}\label{ass.constant.ratio}\\
The ratio between the inductance and resistance of each transmission line is constant, i.e., there exist $\rho>0$ so that
\begin{align*}
\frac{\ell_{jk}}{r_{jk}} = \rho,\quad\forall \, (j,k) \in\mc E. 
\end{align*}
\end{assumption}
{Assumption~\ref{ass.constant.ratio} is usually satisfied for transmission lines at the same voltage level. It is not, however, satisfied across different voltage levels.  While this assumption is required to derive the theoretical results in the remainder of the paper, we observe in simulation that stability is not compromised if the lines have heterogeneous ${\ell}/{r}$ ratio.}   Let us define the quantity $\kappa  \coloneqq \tan^{-1}(\rho\,\omega_0)$. Given the power system's graph $\mc G$, it is easy to verify that the extended Laplacian $\mc L \coloneqq L\otimes I_2$ can be constructed as
\begin{align}\label{eq.laplacian.powernetwork}
 \mc L= \mc R(\kappa) \mc Y,
\end{align}
where {$\mc R(\kappa) = I_N \otimes R(\kappa)$}. The extended Laplacian $\mc L $ is equivalent to the Laplacian of a larger graph that presents two identical connected components, one for the $\alpha$ and one for the $\beta$ coordinate. By using~\eqref{eq.laplacian.powernetwork}, we note that 
\begin{align}\label{eq.laplacian.powernetwork.two}
y = \mc R(\kappa) i_{o} = \mc L \,v.
\end{align}
is a relative output signal (as in~\eqref{eq.oscillators}) that can be obtained from local current measurements.

\subsection{Decentralized control of inverter-based power grids}
Due to abstracting power inverters as controllable voltage sources, we can prescribe any voltage dynamics that relies only on local measurement to the terminal voltage. We choose the following controller structure depicted in Figure~\ref{fig.inverter.controlled}:
\[
\frac{\diff}{\diff t} v_k   = u_k(v_k,y_k), \quad k\in\mc V,\\
\]
where {each component $y_k \in \mathbb{R}^2$ of the network variable $y \in \mathbb{R}^{2N}$ defined in~\eqref{eq.laplacian.powernetwork.two} is a function of the local output current $i_{o,k} \in \mathbb{R}^2$ and can be obtained locally by each controller.} The problem of controlling power inverters falls in the problem class of controlling coupled oscillators of the form~\eqref{eq.oscillators}. Next, we argue that the control objectives~\eqref{eq.obj.theta}~-~\eqref{eq.obj.norm} also have a natural interpretation in terms of power systems. To this end, we introduce the following definition of \emph{instantaneous} active and reactive power (see \cite{AK+83}).

\begin{definition}{\bf(Instantaneous Power)}\label{def:IP}\\
 Given the AC voltage $v_k$ at a node $k \in \mc V$ and the AC current $i_{o,k}$ flowing out of the node, we define the corresponding \emph{instantaneous active power} $p_k \coloneqq v^\top_{k} i_{o,k} \in \mathbb{R}$ and \emph{instantaneous reactive power} $q_k  \coloneqq v^\top_{k}J i_{o,k} \in \mathbb{R}$ flowing out of the node.
\end{definition}

{In the context of controlling an inverter-based grid, the control objective is to steer} the voltage vectors to a synchronous harmonic solution where each inverter operates at the prescribed voltage magnitude and injects the prescribed power into the grid. That is, given active and reactive power set-points $p_k^\star$ and $q_k^\star$ we want the system to converge to a steady state solution in $\mc F_{f,\omega_0} \cap \mc A \cap \mc P $, where $\mc F_{f,\omega_0}$ and $\mc A$ are defined in~\eqref{eq.obj.theta} and~\eqref{eq.obj.norm}  and 
\[
\mc P\coloneqq\left\{ v \in \R^{2N}  \left\vert\;  i_o =\mc Y v,
\begin{bmatrix} v_{k}^\top  \\  v_{k}^\top J \end{bmatrix} i_{o,k} = \begin{bmatrix} p^\star_k \\ q^\star_k\end{bmatrix},\forall\,k\in\mc V
\right.\right\}.
\]

Let us now give a formal definition of a feasible power-flow solution.
\begin{definition}\label{def.fesasible}
A set of power and voltage set-points $\{(p^\star_k,q^\star_k,{v^\star_k})\}_{i\in\mc V}$ is \emph{feasible} if there exist angle set-points $\theta^{{\star}}_{jk}$ such that the power-flow equations 

\begin{align}\label{eq.powers}
\begin{split}
p^\star_k &= \!\!\!\sum_{j:(j,k)\in\mc E}  \!\!\frac{ {v_k^{\star2}} {r_{jk}}  - {v_k^{\star}}{v_j^{\star}}( {r_{jk}}\cos(\theta_{jk}^\star )\!+{ \omega_0\ell_{jk}} \sin(\theta_{jk}^\star )) }{r_{jk}^2 + \omega_0^2 \ell_{jk}^2}{,} \\
q^\star_k &= \!\!\!\sum_{j:(j,k)\in\mc E} \!\!\! \frac{ {v_k^{\star2}} { \omega_0\ell_{jk}} \!-\! {v_k^{\star}}{v_j^{\star}}( { \omega_0\ell_{jk}}\cos(\theta_{jk}^\star )\!-{ r_{jk}} \sin(\theta_{jk}^\star ))}{r_{jk}^2 + \omega_0^2 \ell_{jk}^2}
\end{split}
\end{align}
are satisfied.
\end{definition}

Equation~\eqref{eq.powers} is derived by expanding the standard power-flow equations~\cite{kundur1994power} in terms of the line resistances and inductances and ensures the existence of a corresponding steady-state of an inverter based power system \cite{GD17}.

For fixed voltage magnitudes, a feasible solution to the power flow equations according to Definition~\ref{def.fesasible}, can be parametrized by relative phase shifts $\{\theta_{{k}1}^\star\}_{{k}=2}^N$. The active and reactive power set-points are related to the angle and voltage set-points by the power-flow equations~\eqref{eq.powers}.  Therefore, by choosing $\{\theta_{{k}1}^\star\}_{{k}=2}^N$ consistently with $p_k^\star$, $q_k^\star$, and ${v^\star_k}$ according to~\eqref{eq.powers}, the control specifications \eqref{eq.obj.theta}~-~\eqref{eq.obj.norm} have a natural interpretation in terms of system-wide objectives: If the set $\mc T_{f,\omega_0}$ is asymptotically stable, all state variables in the power network converge to the pre-specified solution of the power-flow equations and, at steady-state, evolve as sinusoidal signals with synchronous frequency $\omega_0$. 

\subsection{Decentralized implementation for networks of inverters}\label{sec.power.systems.implementation}
There is still one obstacle to a decentralized implementation for the power-system application: the matrices $K_k$ in the control law~\eqref{eq.controller} defined in~\eqref{eq.Kk.def} depend on all the angle set-points $\theta_{jk}^\star$ and all the line parameters $\ell_{jk}$ and $r_{jk}$ through the coefficients $w_{jk}$. For power-systems applications this is undesirable as set-points to the inverters are not given in terms of relative angles but powers (active and reactive) and voltage magnitudes. Furthermore, the line parameters of the whole network are, in general, not known by the inverter operator. In this section, we show that the matrices $K_k$, necessary to compute the controller~\eqref{eq.controller}, can be constructed only knowing the power and voltage set-points, and the inductance-resistance ratio $\rho$.

By manipulating the power-flow equations~\eqref{eq.powers}, we obtain a description of $K_k$ that is independent of both the angle set-points and the line parameters and is therefore realistically implementable. 
\begin{proposition}{\bf(Implementable controller)}\label{propo.implementable}\\
Given feasible set-points $\{(p^\star_k,q^\star_k,{v^\star_k})\}_{i\in\mc V}$ according to Definition~\ref{def.fesasible}, the matrices 
$
K_k  \coloneqq  \sum_{j:(j,k)\in\mc E}  w_{jk}{(I_2-\frac{v^\star_j}{v^\star_k}R(\theta_{jk}^\star))}
$
can be constructed as
\begin{align}\label{eq.K.in.PQ}
K_k =  \frac{1}{{v^{\star2}_k}}R(\kappa) \begin{bmatrix}
p_k^\star & q_k^\star\\
-q_k^\star & p_k^\star
\end{bmatrix}.
\end{align}
\end{proposition}
The proof is given in the Appendix. Considering Proposition~\ref{propo.implementable} we can implement the controller~\eqref{eq.controller} using only $p_k^\star, q_k^\star$ and $v^\star_k$. Note that $p_k^\star, q_k^\star$ and $v^\star_k$ need to be consistent across the network, i.e., solve the power-flow equations. \red{We assume that, similar to conventional thermal generation, the set-points are updated relatively infrequently by a centralized decision maker to account for changes in demand and generation. Therefore, on the time-scales of interest for stability analysis we consider the set-points to be constant. This assumption is also universally made in the analysis of established control algorithms for inverters such as droop control \cite{JMG-MC-TLL-PLC:13a}.}
\red{
\begin{remark}
In large scale power systems the set-points $p_k^\star, q_k^\star$ and $v^\star_k$ can be computed by solving a so-called optimal power flow problem that accounts for economic optimality, security, and network congestion \cite{CR+11}. In practice the set-points will only approximately solve the power-flow equations \eqref{eq.powers} due parameter uncertainty and changes in load and generation. In future systems with large share of (fluctuating) renewable generation more frequent updates to the set-points may be required. Several approaches based on distributed optimization and online optimization have been proposed to tackle this problem \cite{DZG13,DS18,TSL17,HZ+17}.
\end{remark}
}
\red{Robustness against the uncertainties in the set-point is implied by the so-called droop behavior of synchronous machines and control algorithms for inverters. In the next section, we discuss the droop properties of the controller \eqref{eq.controller}.}

\subsection{Droop-like behavior}
Power droop is a control law for power inverters where frequency and voltage are adjusted by each inverter proportionally to the local \red{imbalance of active and reactive power}.  It is well known in the literature~\cite{simpson2013synchronization,dorfler2013synchronization,dorfler2012synchronization,MS-FD-BJ-SD:14b} that inverters controlled in this way exhibit dynamics that can be modeled as a network of coupled Kuramoto oscillators~\eqref{eq.kuramoto}. In Section~\ref{sec.kuramoto} we show that, by restricting the closed loop dynamics~\eqref{eq.desired.cl.magnitude} to the circle by means of a projection operation, we recover precisely the nonlinear Kuramoto model. It is then not surprising that, close to the nominal operating point and purely inductive lines (i.e., $\kappa = \frac{\pi}{2}$), the controller~\eqref{eq.controller} resembles the classic droop characteristics.  Let us define 
$\nu_k \coloneqq \norm{ v_k}$ and recall the definition of {$\theta_k \in \mathbb{S}^1$ as the angle of the voltage $v_k$, i.e., $\theta_k$ satisfies $v_k = R(\theta_k) ({\norm{v_k}},0)$}.
\begin{proposition}\label{propo.droop.like} {\bf (Polar coordinates)}
If $r_{jk}=0,~\forall \,(j,k)\in\mc E$, the closed loop~\eqref{eq.desired.cl.magnitude} in polar coordinates is given by
\begin{align}
\begin{split}\label{eq.droop.like}
\dot \nu_k = &~ \eta\left(\frac{q^\star_k}{{v^{\star2}_k}} - \frac{q_k}{\nu_{k}^2}\right)  \nu_k + \frac{\alpha}{{v^\star_k}}({v^\star_k}-\nu_k)\nu_k\\
\dot \theta_k = &~  \omega_0+\eta \left( \frac{p_k^\star}{{v^{\star2}_k}} - \frac{p_k}{\nu_{k}^2} \right)
\end{split}
\end{align}
where $p_k = {v^\top_k i_{o,k}}$ and $q_k = {v^\top_k J i_{o,k}}$ are the active and reactive powers injected by the $k^\text{th}$ inverter.
where $p_k = {v^\top_k i_{o,k}}$ and $q_k = {v^\top_k J i_{o,k}}$ are the active and reactive powers injected by the $k^\text{th}$ inverter.
\end{proposition}
The proof is provided in the Appendix. Equation~\eqref{eq.droop.like} shows that the consensus-inspired controller~\eqref{eq.controller}, similarly to virtual oscillator control~\cite{MS-FD-BJ-SD:14b}, resembles a droop controller around the normal operating point ($\nu_k\approx {v^\star_k}=1$). \red{Moreover, the synchronization and active power-sharing properties of droop control that imply robustness against inconsistent active power set-points can only be rigorously established for inductive lines, constant voltage magnitude (i.e., $\dot \nu_k = 0$), and near the nominal operating point \cite{simpson2013synchronization}. Under these restrictive assumptions \eqref{eq.droop.like} is equivalent to standard active power droop control with a re-scaled power set-point and has the same properties.}
However,~\eqref{eq.droop.like} differs greatly from a standard droop dynamics during large transients. In contrast to droop control, this allows us to prove almost global asymptotic stability of the grid with respect to consistent set-points from almost all initial conditions (see Theorem~\ref{thm.ags} and Proposition~\ref{propo.zero.unstable}).

\red{
Moreover, as shown in \cite[{Prop. 3}]{colombino2017global}, for purely resistive lines (i.e., $\ell_{jk}=0,~\forall \,(j,k)\in\mc E$) the closed-loop system \eqref{eq.desired.cl.magnitude.two} rewritten in polar coordinates is similar to the averaged equations of the angle and voltage dynamics of virtual oscillator control \cite{johnson2016synthesizing}. However, unlike virtual oscillator control, the controller proposed in this work is fully dispatchable and can drive the system to a prescribed solution of the power-flow equations. In contrast, virtual oscillator control synchronizes a network of inverters to zero phase difference and identical magnitudes and is robust with respect to unknown (linear) passive loads \cite{LABT-JPH-JM:13}. As a consequence virtual oscillator control delivers power to the loads, but the voltage magnitude and the power sharing between the inverters cannot be controlled.}

\red{Given that the control \eqref{eq.controller} closely resembles the established methods in the inductive and resistive limit,} we expect that the robustness properties of the proposed controller with respect to inconsistency in the set-points are similar to established methods that rely on virtual Li\'enard-type oscillators and so-called droop characteristics, i.e., a trade-off between power imbalance and deviations of frequency and voltage. \red{Preliminary simulation results that confirm the droop characteristics of the controller \eqref{eq.controller} can be found in \cite{SC+19}. A rigorous analysis is an interesting topic of future work.}

\section {Analysis of the closed-loop dynamics}\label{sec.proofs}
In Section~\ref{sec.consensus.control} we introduced the desired closed-loop dynamics and showed that they can be obtained with a decentralized implementation. In order to state the main result of the paper, we begin by defining the precise notion of stability with respect to a set. Next, we provide some preliminary results that culminate with Theorem~\ref{thm.ags} that shows that under {an insightful and tractable stability condition} on the angle set-points and grid parameters, the closed-loop dynamics~\eqref{eq.desired.cl.magnitude} are {almost globally asymptotically stable} from almost all initial conditions.

\subsection{Basic definitions}
We begin by defining stability and almost global attractivity of a dynamic system with respect to a set.
\begin{definition}\label{def:stability}
 A dynamic system $\ddt x = f(x)$ is called stable with respect to a set $\mc C$ if for every $\epsilon \in \mathbb{R}_{>0}$ there exists $\delta=\delta(\epsilon)\in \mathbb{R}_{>0}$ such that
  \begin{align}
   \norm{x_0}_{\mc C} < \delta \implies \norm{\varphi(t,x_0)}_{\mc C} < \epsilon, \qquad \forall t \geq 0.
   \end{align}
\end{definition}
\begin{definition}
 A dynamic system $\ddt x = f(x)$ is called almost globally attractive with respect to a set $\mc C$ if for every pair $\epsilon \in \mathbb{R}_{>0}$ and $\delta \in \mathbb{R}_{>0}$ there exists $T=T(\epsilon,\delta) \in \mathbb{R}_{>0}$ such that for all $x_0 \notin \mc Z$
  \begin{align}\label{eq:def:attr}
     \norm{x_0}_{\mc C} < \delta \implies  \norm{\varphi(t,x_0)}_{\mc C} < \epsilon, \quad \forall t \geq T,
   \end{align}
       and $\mc Z$ has zero Lebesgue measure. If, in addition, $\mc Z = \emptyset$, the system is called globally attractive.
\end{definition}
Almost global asymptotic stability with respect to a set $\mc C$ is defined by combining stability and almost global attractivity (see \cite{A04}).
\begin{definition}\label{def:ags}
A dynamic system $\ddt x = f(x)$ is called (almost) globally asymptotically stable with respect to a set $\mc C$ if it is stable with respect to $\mc C$ and (almost) globally attractive with respect to $\mc C$. A dynamic system $\ddt x = f(x)$ is called globally exponentially stable with respect to a set $\mc C$ if it is globally stable with respect to $\mc C$ and exponentially converges to $\mc C$.
\end{definition}
As is customary for autonomous systems, these definitions imply uniform stability and uniform attractivity. Observe, that in contrast to the definition in \cite{A04}, we do not require the set $\mc C$ to be compact. Finally, the standard notion of global asymptotic stability \cite{H67} can be recovered by letting $\mc C = \{0\}$ and $\mc Z = \emptyset$.

\subsection{Stability condition}
In the following, we provide a condition on the angle set-points and control gains that is supported by physical intuition and guarantees that~\eqref{eq.desired.cl.magnitude} is almost globally asymptotically stable with respect to the control objectives. 
\begin{condition}\label{cond.small.angle}{\bf(Stability condition)}\\
The graph $\mc G$ is connected and $\theta^\star_{k1} \in {[}0,\frac{\pi}{2}{]}$ holds for all $k \in \mc V \setminus \{1\}$. Moreover, the angle set-points $\{{\theta_{k1}^\star}\}_{k=1}^N$ as well as the design parameters $\alpha>0$ and $\eta>0$ satisfy
\begin{align}\label{eq.condtion.js}
\max_k \sum_{j=1}^N w_{jk}\left|1-\frac{v^\star_j}{v^\star_k}\cos(\theta^\star_{jk})\right| + \frac{\alpha}{\eta} < \frac{1}{2}\frac{v^{\star2}_{\min}}{v^{\star2}_{\max}}\lambda_2( L),
\end{align}
where $v^{\star}_{\min} \coloneqq \min_{k \in \mc V} v^\star_k$ and $v^{\star}_{\max} \coloneqq \max_{k \in \mc V} v^\star_k$ are respectively the smallest and largest magnitude set-points and $\lambda_2(L)$ is the second smallest eigenvalue of the graph Laplacian $L$.
\end{condition}
If the graph is connected (thus $\lambda_2(L)>0$), the inequality \eqref{eq.condtion.js} can always be satisfied by a proper choice of control gains and small enough angle set-points $\{{\theta_{k1}^\star}\}_{k=1}^N$.

\begin{remark}{\bf{(Interpretation of the stability condition)}}\\
Condition~\eqref{eq.condtion.js} can be interpreted as a heterogeneity-connectivity trade-off for the desired steady state, and it is closely related to the stability conditions for Kuramoto oscillators presented in~\cite{dorfler2013synchronization}. 

In the power systems context presented in Section~\ref{sec.inverter.model}, under the assumption that the transmission lines are purely inductive, i.e., $r_{jk}=0$, the stability condition~\eqref{eq.condtion.js} becomes 
\begin{align*}\label{eq.condition.ps}
 \max_k \sum_{j=1}^N\frac{1}{{v^{\star2}_k}}\left| q^\star_{jk}\right|+\frac{\alpha}{\eta} <\frac{1}{2}\frac{v^{\star2}_{\min}}{v^{\star2}_{\max}}\lambda_2( L),
\end{align*}
where $q^\star_{jk}$ are the reactive branch powers of the steady state solution of the power flow {(see \eqref{eq.powers.ref})}. In the literature of virtual oscillator control for power inverters~{\cite{johnson2014synchronization,johnson2016synthesizing,MS-FD-BJ-SD:14b,sinha2016synchronization,LABT-JPH-JM:13}} it is shown that, for connected graphs, {global synchronization can be obtained. However, VOC cannot be dispatched to track pre-specified power and voltage setpoints.} Condition~\eqref{eq.condition.ps} (together with Theorem~\ref{thm.ags} in Section~\ref{sec.main.result}) ensures that {a pre-specified solution of the power-flow equations} (with nonzero power flows) is almost globally stable.
\end{remark}

\subsection{Main result: Almost global {asymptotic stability}}\label{sec.main.result}
This section is devoted to showing that the trajectories of the closed-loop system~\eqref{eq.desired.cl.magnitude} indeed converge to the invariant set $\mc T_{f,\omega_0}${, i.e.,} the controller~\eqref{eq.controller} drives the system to a synchronous steady-state with the correct angle set-points (power flows for the powers systems application) and amplitudes.

\begin{theorem}{\bf (Almost global asymptotic stability)}\label{thm.ags}\\
Under Condition~\ref{cond.small.angle}, the dynamical system~\eqref{eq.desired.cl.magnitude} is almost globally asymptotically stable with respect to the set $\mc T_{f,\omega_0}$, i.e., with respect to the target sets specifying the control objectives~\eqref{eq.obj.theta},~\eqref{eq.obj.sync}, and~\eqref{eq.obj.norm}.
\end{theorem}

The remainder of this section is devoted to building the necessary auxiliary results that allow us to prove Theorem~\ref{thm.ags}. In view of Definition \ref{def:ags}, in order to prove Theorem~\ref{thm.ags}, we need to show that the target set is both (almost) globally attractive and stable. 

We begin by considering both the closed-loop system and the control objectives in a rotating reference frame so that equilibria in the new frame correspond to synchronous solutions in the static frame. In the rotating frame, the target set reduces to the intersection of the set $\mc S$ of states with correct relative magnitudes {and phase} defined in~\eqref{eq.obj.theta} and the set $\mc A$ of states with correct magnitude defined in~\eqref{eq.obj.norm}. 

\begin{figure}[htbp]
\begin{center} 
\includegraphics[width=\columnwidth]{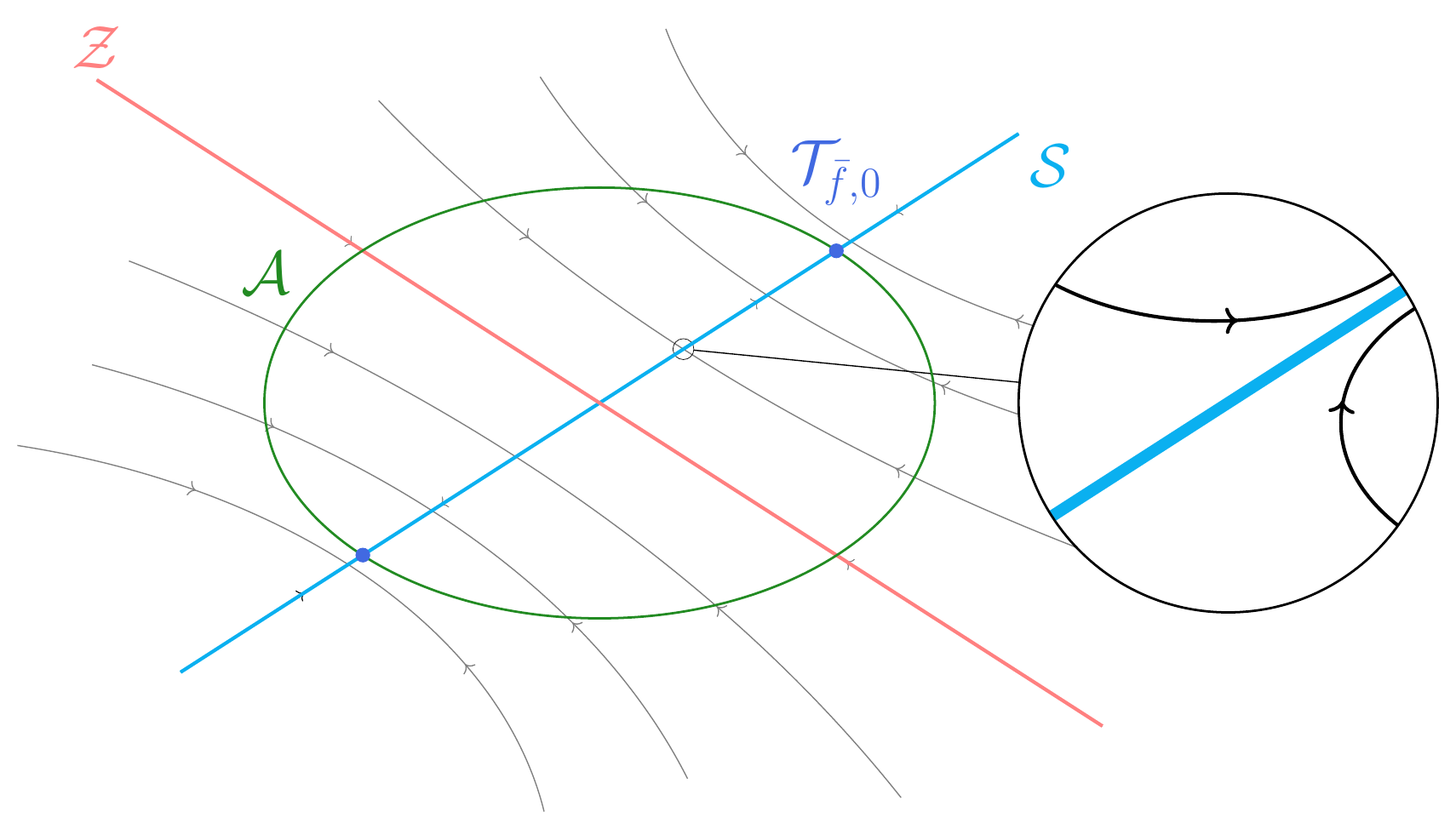}
\end{center}
\caption{The trajectories $\bar v(t)$ approach the subspace $\mc S$ exponentially fast. Once the distance $\norm{\bar v}_{\mc S}$ of $\bar v$ from $\mc S$ is sufficiently small, provided that $v(t)\in\R^n\backslash\mc Z$, where $\mc Z$ is a set of zero Lebesgue measure, the magnitude controller steers the trajectories towards $\mc A$. As a result all trajectories except those starting on $\mc Z$ will approach $\mc S \cap \mc A =\mc T_{\bar{f},0}$. \label{fig.convergence}}
\label{fig.attractivity}
\end{figure}

As a first step to proving (almost) global attractivity Theorem~\ref{th.convergence.to.S} establishes that the set $\mc S$ is globally exponentially stable. In Proposition~\ref{propo.wzeroinputconv} we prove that $\mc A$ is asymptotically stable  when the dynamics are restricted to $\mc S$. Proposition~\ref{propo.zero.unstable} shows that the origin, which is the only equilibrium of the system that lies outside the target set, has a zero-measure region of attraction. {In} Theorem~\ref{thm.sync} we merge these results using a continuity argument and we prove that, as illustrated in Figure~\ref{fig.attractivity}, almost all trajectories approach $\mc S$ exponentially and eventually approach $\mc A$. In addition the proof of Theorem~\ref{thm.sync} requires Proposition~\ref{propo.wbound}, which guarantees that the trajectories remain bounded at all times.

Finally, with Proposition~\ref{propo.invset} we construct a family of nested invariant sets that will be used in Theorem~\ref{thm.stab} to establish stability of the target set.

\subsubsection{Proof of attractivity}
We begin by performing a change of variables for the closed loop~\eqref{eq.desired.cl.magnitude} to a reference frame rotating at speed $\omega_0$:
\begin{align*}
\bar v = \diag( \{R(\omega_0 t)\}_{k=1}^N ) v.
\end{align*}
Since $\diag( \{R(\omega_0 t)\}_{k=1}^N )$ commutes with both $\mc K -\mc L$ and $\Phi(v)$, the dynamics~\eqref{eq.desired.cl.magnitude} in the new coordinates become
\begin{align}\label{eq.desired.cl.magnitude.rotframe}
\frac{\diff}{\diff t} \bar v = \eta (\mc K - \mc L) \bar v + \alpha \; \Phi(\bar v)\, \bar v \eqqcolon \bar{f}(\bar{v}).
\end{align}
\begin{proposition}{\bf{(Control objectives in the rotating frame)}}\label{propo.equivalent.rotating}\\
The following statements are equivalent:
\begin{enumerate}
\item The dynamics \eqref{eq.desired.cl.magnitude} are almost globally asymptotically stable with respect to $\mc T_{f,\omega_0}$.
\item The dynamics \eqref{eq.desired.cl.magnitude.rotframe} are almost globally asymptotically stable with respect to $\mc T_{\bar{f},0}$.
\end{enumerate}
\end{proposition}
The Proposition directly follows by noting that the transformation into the rotating frame is a diffeomorphism. Proposition~\ref{propo.equivalent.rotating} allows us to analyze the closed loop in the rotating coordinate frame. 

Let us begin by showing that $\mc S$ is exponentially stable for the closed-loop dynamics~\eqref{eq.desired.cl.magnitude.rotframe}. To this end, we define the matrix whose span is the set $\mc S$ as
\begin{align*}
S\coloneqq \frac{1}{\sqrt{\sum_k v_k^{\star2}}}\begin{bmatrix}
v^\star_1 I_2\\
v^\star_2 R(\theta_{21}^\star)\\
\vdots\\
v^\star_N R(\theta_{N1}^\star)
\end{bmatrix},
\end{align*}
and we define the matrix $P$ as the projector onto the subspace orthogonal to $\mc S$ as $P\coloneqq I_{2N}-SS^\top$ {and the Lyapunov function
\begin{align}\label{eq.lyap.def}
 V(\bar{v}) =  \bar{v}^\top P \bar{v}  =  \norm{\bar{v}}^2_{\mc S}.
\end{align}
}
\begin{proposition} \textbf{(Bound on Lyapunov derivative)}\label{propo.bounded}\\
The time derivative of $V(\bar{v})$ along the trajectories of~\eqref{eq.desired.cl.magnitude.rotframe} is bounded by
\begin{align}\label{eq.lyap.decr}
\!\!\frac{\diff}{\diff t}{V}(\bar{v}) \!\le\! &~ \bar{v}^\top\! \left[\eta \left(P ( \mc K - \mc L ) +(\mc K - \mc L )^\top P \right)+ 2\alpha P\right] \bar{v}.
\end{align}
\end{proposition}
\begin{IEEEproof}
Taking the time derivative of $V(\bar{v})$ along the trajectories of~\eqref{eq.desired.cl.magnitude.rotframe} we obtain
\begin{align*}\label{eq.lyap.dec.mag}
\frac{\diff}{\diff t}{V}(\bar{v}) = &~\eta \bar{v}^\top \left[P ( \mc K - \mc L) +(\mc K - \mc L )^\top P\right] \bar{v} \\
+&~ \alpha\; \bar{v}^\top \left[P\; \Phi(\bar{v})+ \Phi(\bar{v})P\right]\bar{v}.
\end{align*}
In order to prove the claim we show that 
\begin{align}\label{eq.claim1}
\bar{v}^\top P \Phi(\bar{v}) \bar{v} \le \bar{v}^\top P \bar{v},\quad \forall\, \bar{v}\in\R^{n}.
\end{align}
Since
\[
\bar{v}^\top P \Phi(\bar{v}) \bar{v}  =  \bar{v}^\top P \bar{v} -  \bar{v}^\top  \diag\left(\left\{\frac{1}{v_k^\star}\norm{\bar{v}_k} I_2\right\}_{k=1}^N\right) P\,\bar{v},
\]
the inequality~\eqref{eq.claim1} is equivalent to  
\[
  \bar{v}^\top  \diag\left(\left\{\frac{1}{v_k^\star}\norm{\bar{v}_k} I_2\right\}_{k=1}^N\right) P\,\bar{v} \geq 0, \quad \forall \bar{v}\in\R^{n}.
\]
{
Next, we use $\theta_{kj} = \theta_k - \theta_j \in \mathbb{S}^1$ to denote the relative angle between $v_j$ and $v_k$ such that $\norm{v_k} v_j = \norm{v_j} R(\theta_{kj}) v_k$ holds. Given the particular form of $P$ we can write} 
\begin{align*}
& \sum_n v_n^{\star2}~ \bar{v}^\top  \diag\left(\left\{\frac{1}{v_k^\star}\norm{\bar{v}_k} I_2\right\}_{k=1}^N\right) P\,\bar{v}\\
& = \sum_{k=1}^N \bigg(\frac{\sum_n v_n^{\star2}}{{v^\star_k}}\|\bar{v}_k\|^3 -  \sum_{j=1}^N v^\star_j\|\bar{v}_k\|^2\norm{\bar{v}_j}\cos(\theta_{kj} - \theta^\star_{kj})\bigg) \\
&\ge \sum_{k=1}^N \left(\frac{\sum_n v_n^{\star2}}{{v^\star_k}}\|\bar{v}_k\|^3 -  \sum_{j=1}^N v^\star_j\|\bar{v}_k\|^2\norm{\bar{v}_j}\right)\\
&= \sum_{k=1}^N \left(\frac{\sum_n v_n^{\star2}}{{v^\star_k}}\|\bar{v}_k\|^2 \left( {\|\bar{v}_k\|} -{v_k^\star}\, \tilde v \right)\right),
\end{align*}
where
$
\tilde v\coloneqq \frac{1}{\sum_n v_n^{\star2}}\sum_{j=1}^N v^\star_j\norm{\bar{v}_j}\ge 0.
$
Next, let us define 
\[
 A \coloneqq \sum\nolimits_{k=1}^N \left({\sum_n v_n^{\star2}}\|\bar{v}_k\| \left( {\|\bar{v}_k\|} -{v_k^\star}\, \tilde v \right)\right),
\]
and
\[
 B \coloneqq \sum\nolimits_{k=1}^N \frac{\sum_n v_n^{\star2}\|\bar v_k\|}{{v^\star_k}}  \left( {\|\bar{v}_k\|} -{v_k^\star}\, \tilde v \right)^2.
\]
Letting 
$
s \coloneqq \begin{bmatrix}
v^\star_1&
\cdots&
v^\star_N
\end{bmatrix}^\top\!\!,
$
it follows that 
\[
A =  [\norm{\bar v_1},...,\norm{\bar v_n}] \left(\sum_n v_n^{\star2} I_N - ss^\top\right) 
\begin{bmatrix}
\norm{\bar v_1}\\
\vdots\\
\norm{\bar v_N}
\end{bmatrix}\ge 0.
\]
Moreover, $B\ge0$ holds because $B$ it is a sum of nonnegative quantities. One can easily verify that
\begin{align*}
& \sum\nolimits_{k=1}^N \left(\frac{\sum_n v_n^{\star2}}{{v^\star_k}}\|\bar{v}_k\|^2 \left( {\|\bar{v}_k\|} -{v_k^\star}\, \tilde v \right)\right) = \tilde v A + B,
\end{align*}
and it immediately follows that the inequalities \eqref{eq.claim1} and \eqref{eq.lyap.decr} hold.
\end{IEEEproof}
Next, we show that Condition \ref{cond.small.angle} ensures that $V(\bar{v})$ is deceasing for all $\bar{v} \notin \mc S$. 
\begin{proposition}\textbf{(Lyapunov decrease on $\mathbf{\mc S^\perp}$)}\label{propo.decreasing}\\
Under Condition~\ref{cond.small.angle} it holds that
\begin{align}\label{eq.ass.reformulation}
v^\top \left[\eta[(\mc K-\mc L)^\top P+P(\mc K-\mc L)] + 2\alpha P \right] v<0 \quad \forall v\in\mc S^\perp.
\end{align}
\end{proposition}
\begin{IEEEproof}
We start by noting that, since $(\mc K-\mc L){v}=0$ for all $v\in\mc S$ and $P$ is the projector onto $\mc S ^\perp$, it holds that $(\mc K-\mc L) = (\mc K-\mc L)P.$
Therefore,~\eqref{eq.ass.reformulation} is equivalent to, 
\begin{align}
 & v^\top[\eta  P(\mc K-\mc L) P+\alpha P] v <0\quad \forall v\in\R^n \nonumber \\ 
\iff   &v^\top [\eta P\mc KP +\alpha P] v<v^\top \eta P\mc LPv \quad \forall v\in\R^n \label{eq.ineq.condition}
\end{align}
For all $v\in\R^n$, we can bound the left hand side of~\eqref{eq.ineq.condition} by noting that 
\begin{align}
v^\top P\mc KP v&\le \frac{1}{2}\|\mc K+\mc K^\top\|\|v\|_{\mc S}^2\nonumber\\
& = \|v\|_{\mc S}^2\max_k \sum_{j=1}^N w_{jk}\left|1-\frac{v^\star_j}{v^\star_k}\cos(\theta^\star_{jk})\right|,  \label{eq.bound0} 
\end{align}
where~\eqref{eq.bound0} follows from the fact that $\frac{1}{2}(\mc K +\mc K^\top)=\diag\left (\left\{  \sum_{j=1}^N w_{jk}\left(1-\frac{v^\star_j}{v^\star_k}\cos(\theta^\star_{jk})\right) \right\}_{k=1}^N\right)$. 

In view of~\eqref{eq.bound0}, the left-hand side of~\eqref{eq.ineq.condition} can be bounded as follows
\begin{align}
&v^\top [\eta P\mc KP +\alpha P] v \notag \\ &\quad \le  \|v\|_{\mc S}^2\left[\eta \max_k \sum_{j=1}^N w_{jk}\left|1-\frac{v^\star_j}{v^\star_k}\cos(\theta^\star_{jk})\right|+\alpha\right]. \label{eq.bound3}
\end{align}
Before bounding the right-hand side of~\eqref{eq.ineq.condition}, we define the matrix $\mathbb{I}\in\R^{2N\times 2}$ as
\[
\mathbb I = \begin{bmatrix}I_2 & \cdots & I_2\end{bmatrix}^\top,
\]
and the projector onto the orthogonal subspace of the span of the matrix $\mathbb I$ as
\[
P_{\mathbb I^\perp}\coloneqq I_{2N}-\frac{1}{N}\mathbb{I\,I}^\top.
\]
Note that both $P$ and $P_{\mathbb I^\perp}$ are projector matrices, therefore $P^2=P$ and $P^2_{\mathbb I^\perp}= P_{\mathbb I^\perp}$. Furthermore, since $\mc Lv= 0$ for all $v\in\mathrm{range}(\mathbb I)$, it holds that $\mc LP_{\mathbb I^\perp} = \mc L$. We can now bound the right hand side of~\eqref{eq.ineq.condition} as 
\begin{align}
&v^\top P\mc LPv   =  v^\top PP_{\mathbb I^\perp}\mc LP_{\mathbb I^\perp}Pv \nonumber \\
  & \ge \lambda_2( L)\|P_{\mathbb I^\perp}Pv\|^2\nonumber\\
  & = \lambda_2( L) \left(v^\top P v -\frac{1}{N}v^\top P \mathbb{I\,I}^\top P v\right)\nonumber\\
    & = \lambda_2( L) \left(v^\top P v -\frac{1}{N}v^\top P^2 \mathbb{I\,I}^\top P^2 v\right)\nonumber\\
  & \ge \lambda_2( L)\|v\|_{\mc S}^2\left(1-\frac{1}{N}\|P\mathbb{I}\|^2\right)\nonumber\\
  & = \lambda_2( L)\|v\|_{\mc S}^2\left(1-\frac{1}{N}\|\mathbb I^\top P \mathbb I\|\right)\nonumber\\
  & = \lambda_2( L)\|v\|_{\mc S}^2\nonumber\\
  &\left(1-\frac{1}{N{\sum_n v^{\star2}_n}}\left\|N{\sum_n v^{\star2}_n}-\sum_{k=1}^N\sum_{j=1}^Nv^\star_j{v^\star_k}R(\theta^\star_{jk})\right\|\right).\label{eq.norm.bound}
 \end{align} 
 Since $R(\theta^\star_{jk})=R(\theta^\star_{kj})^\top$, the off-diagonal terms of the sum 
 $
  \sum_{k=1}^N\sum_{j=1}^Nv^\star_j{v^\star_k}R(\theta^\star_{jk})
 $
 in~\eqref{eq.norm.bound} cancel and it follows that
    \begin{align}
    &N{\sum_n v^{\star2}_n}-\sum_{k=1}^N\sum_{j=1}^Nv^\star_j{v^\star_k}R(\theta^\star_{jk}) \nonumber \\
    & =\left(N{\sum_n v^{\star2}_n}-\sum_{k=1}^N\sum_{j=1}^Nv^\star_j{v^\star_k}\cos(\theta^\star_{jk})\right) I_2,\nonumber \\
    & =\left( \sum_{k=1}^N\sum_{j=1}^N\frac{v^{\star2}_k+v^{\star2}_j}{2} - v^\star_j{v^\star_k}\cos(\theta^\star_{jk}) \right) I_2 \succeq 0, \label{eq.nonnegativity.bound}
\end{align}
where the last inequality follows from the fact that $v^{\star2}_k+v^{\star2}_j\ge 2 v^{\star}_k v^{\star}_j$. Using~\eqref{eq.nonnegativity.bound}, we can express the bound in~\eqref{eq.norm.bound} as 
  \begin{align}
 v^\top P\mc LPv       & \ge \lambda_2( L)\|v\|_{\mc S}^2\frac{1}{N{\sum_n v^{\star2}_n}}  \sum_{k=1}^N\sum_{j=1}^N v^\star_j{v^\star_k} \cos(\theta^\star_{jk}) \nonumber \\ 
  & \ge \frac{1}{2}\frac{v^{\star2}_{\min}}{v^{\star2}_{\max}}\lambda_2( L)\|v\|_{\mc S}^2. \label{eq.bound2}
\end{align}
The last inequality in~\eqref{eq.bound2} follows from the fact that 
\[
 \sum_{k=1}^N\sum_{j=1}^Nv^\star_j{v^\star_k}\cos(\theta^\star_{jk})\ge N^2 v^{\star2}_{\min}\frac{1}{2},
\]
which can be proven by induction by showing that the minimizing choice of angles for the sum of cosines is given by choosing $\{\theta^\star_{k1}\}_{k=2}^N\le \frac{\pi}{2}$ such that $\lceil\frac{N}{2}\rceil$ angles equal to 0 and $N-\lceil\frac{N}{2}\rceil$ equal to $\frac{\pi}{2}$.
By substituting the bounds~\eqref{eq.bound3} and~\eqref{eq.bound2} into~\eqref{eq.ineq.condition}, we can conclude that~\eqref{eq.condtion.js} is a sufficient condition for~\eqref{eq.ass.reformulation}.
\end{IEEEproof}

\begin{theorem}\label{th.convergence.to.S}{\bf(Exponential phase stability)}\\
Given angle set-points $\{{\theta_{k1}^\star}\}_{k=1}^N$ and design parameters $\alpha>0$ and $\eta>0$ that satisfy Condition~\ref{cond.small.angle}, $V(\bar{v}) = \bar{v}^\top P\bar{v}$ is a Lyapunov function for the dynamical system~\eqref{eq.desired.cl.magnitude.rotframe} with respect to $\mc S$. The closed-loop dynamical system~\eqref{eq.desired.cl.magnitude.rotframe} is globally exponentially stable with respect to $\mc S$.
\end{theorem}
\begin{IEEEproof}
Using Propositions~\ref{propo.bounded} and~\ref{propo.decreasing} we conclude that
\begin{align*}
\frac{\diff}{\diff t}{V}(\bar{v}) \le &~ \bar{v}^\top \left[\eta \left(P ( \mc K - \mc L ) +(\mc K - \mc L )^\top P \right)+ 2\alpha P\right] \bar{v} \\
 = &~ \bar{v}^\top Q \bar{v} \le 0 \qquad \forall v\in\mc S^\perp.
\end{align*}
Consider a vector $\xi \in \mathbb{R}^{2N}$. Since $P$ is the projector onto $\mc S ^\perp$, it holds that $\xi \in\mc S$ if and only if $V(\xi)=0$ and $P\xi=Q\xi=0$ hold. {Using $V(\bar{v}) =  \norm{\bar{v}}^2_{\mc S}$ it can be shown that}
\begin{align*}
\frac{\diff}{\diff t}{V}(\bar{v}) &\le -\alpha_2 \norm{\bar{v}}^2_{\mc S} \le -\alpha_2  V(\bar{v}), 
\end{align*}
where $\alpha_2 = \sigma_2(Q)$ is the smallest nonzero singular value of $Q$ and considering \eqref{eq.ass.reformulation} it holds that $\alpha_2 >0$. Using standard results from Lyapunov theory~\cite{H67} it is straightforward to show that $V$ converges to $0$ exponentially and that the dynamical system~\eqref{eq.desired.cl.magnitude.rotframe} is globally exponentially stable with respect to $\mc S$. 
\end{IEEEproof}

We have established that all trajectories converge exponentially to $\mc S$. We now show that the origin, {despite} being in $\mc S$ is unstable and reachable only from a set of measure zero.
\begin{proposition}{\bf(Instability of $\boldsymbol{\bar{v}=0}$)}\label{propo.zero.unstable}\\
  Consider dynamics \eqref{eq.desired.cl.magnitude.rotframe}, the equilibrium $\bar{v}=0$, and $\eta>0$ and $\alpha>0$. Moreover, {consider the set $\mc Z$ in which the origin is attractive}:
  \begin{align*}
   \mc Z \coloneqq \{ \bar{v}_0 \in \mathbb{R}^n \vert \lim\nolimits_{t\to\infty} \varphi_{\bar{f}}(t,\bar{v}_0)=0, \quad \bar{v}(0)=\bar{v}_0 \}.
  \end{align*}
  The set $\mc Z$ has zero Lebesgue measure. Moreover, the equilibrium $\bar{v}=\{0\}$ is unstable.  
 \end{proposition}
 The proof of Proposition~\ref{propo.zero.unstable} is provided in the Appendix. {The zero measure set $\mc Z$ set is characterized by the eigenvectors corresponding to the unstable eigenvalues of the Jacobian $\frac{\partial \bar{f}(\bar{v})}{\partial \bar{v}}\vert_{\bar{v}=0}$. The reader is referred to \cite[Prop. 11]{Monzon2006} and \cite[Sec. 5A]{hirsch2006invariant} for further details.}
 Next, we define the Lyapunov-like functions 
\[
W_k(\bar{v}_k) = \Phi_k(\bar{v}_k)^2,
\]
where $\Phi_k(\bar{v}_k)$ is defined in~\eqref{eq.phi}. With the following proposition we show that the trajectories of the closed loop~\eqref{eq.desired.cl.magnitude.rotframe} remain bounded for all times. 
\begin{proposition}{\bf(Boundedness of the trajectories)}\label{propo.wbound}\\
Under Condition~\ref{cond.small.angle}, for every initial condition $\bar{v}_0 \in \mathbb{R}^n$ and every $k \in \mc V$ there exists a constant {$\hat{W}_k \in \mathbb{R}_{>1}$ such that the compact set 
 \begin{align}\label{eq.sets.wk}
  \mc W_k(\bar{v}_0) = \left\{ \bar{v} \in \mathbb{R} \left.\vert W_k(\bar{v}) \leq \hat{W}_k(\bar{v}_0) \right. \right\}.
 \end{align}
 is invariant with respect to the dynamics \eqref{eq.desired.cl.magnitude.rotframe}.}
\end{proposition}
The proof of Proposition~\ref{propo.wbound} is given in the Appendix.
Finally, we establish that the set $\mc T_{\bar f,0}$ is asymptotically stable in $\mc S \setminus \{0\}$.
\begin{proposition}{\bf(Asymptotic stability of $\boldsymbol{\mc T_{\bar f,0}}$ in $\boldsymbol{\mc S \setminus \{0\}}$)}\label{propo.wzeroinputconv}\\
For all initial conditions $\bar{v}_0 \in \mc S \setminus \{0\}$ the dynamics \eqref{eq.desired.cl.magnitude.rotframe} are asymptotically stable with respect to $\mc T_{\bar f,0}$.
\end{proposition}
\begin{IEEEproof}
It follows from Theorem \ref{th.convergence.to.S} that the set $\mc S$ is positively invariant under the dynamics \eqref{eq.desired.cl.magnitude.rotframe}, i.e, $\bar{v}_0 \in \mc S$ implies $\norm{\varphi_{\bar{f}}(t,\bar{v}_0)}_{\mc S}=0$ for all $t \geq 0$. Substituting $\norm{\bar{v}}_{\mc S}=0$ into \eqref{eq.wk.decr} results in
\begin{align*}
  \frac{\diff}{\diff t} W_k &\leq -2 W_k \norm{\bar{v}_k} \leq 0.
\end{align*}
Next, consider the set of points such that $\bar{v}_k=0$ for some $k\in\mc V$:
\begin{align*}
 \mc O = \left\{ \bar{v} \in \mathbb{R}^n \right\vert\left. \exists k \in \mc V: \bar{v}_k = 0 \right\}.
\end{align*}
By construction of $\mc S$, it holds that $\mc O \cap \mc S = \{0\}$. It follows that $\tfrac{\diff}{\diff t}  W_k  < 0$ for all $\bar{v} \in \mc S \setminus \{0\}$ such that $W_k (\bar{v}_k) > 0$. Using $v^{\star2}_k W_k = \norm{\bar{v}_k}^2_{\mc A_k}$, it can be seen that $\varphi(t,\bar{v}_0) \neq 0$ holds for all $\bar{v}_0 \in \mc S \setminus \{0\}$. In other words, $\mc S \setminus \{0\}$ is invariant under the dynamics \eqref{eq.desired.cl.magnitude.rotframe}. Moreover, using $ v^{\star2}_k W_k = \norm{\bar{v}_k}^2_{\mc A_k}$, it follows from standard results from Lyapunov theory that $\tfrac{\diff}{\diff t}  W_k < 0$ for all $W_k \neq 0$ implies asymptotic stability of $\mc A$ in $\mc S \setminus \mc \{0\}$ (cf. Theorem 42.4 and Section 45 in \cite{H67}). 
\end{IEEEproof}

We are now ready for the two major results that will allow us to prove Theorem~\ref{thm.ags}: We will first establish almost global attractivity of the closed-loop dynamics with respect to $\mc T_{\bar{f},0}$, i.e., the trajectories of the closed-loop dynamics in the rotating reference frame~\eqref{eq.desired.cl.magnitude.rotframe} converge to $\mc T_{\bar{f},0}$ from almost every initial condition. Figure~\ref{fig.convergence} illustrates the main idea behind the proof. We exploit the fact that for all initial conditions the system converges exponentially to $\mc S$ driving the {state} vectors to the correct relative {phase} angles. {Using similar arguments as in \cite{SontagCICS03} we show} that, unless the trajectory starts from the {zero-measure set} $\mc Z$, once the system is close enough to $\mc S$, the {states $v_k$ converge to the correct magnitude.}
{
\begin{theorem}{\bf (Almost global attractivity)}\label{thm.sync}\\
 Under Condition~\ref{cond.small.angle}, the dynamical system~\eqref{eq.desired.cl.magnitude.rotframe} is almost globally uniformly attractive with respect to $\mc T_{\bar f,0}$.
 \end{theorem}}
 \begin{IEEEproof}
We prove the theorem by showing that $\mc T_{\bar{f},0}$ is {uniformly} attractive for all $\bar{v}_0 \in \R^n \setminus \mc Z$. In other words, we show that for every pair $\delta \in \R_{>0}$ and $\epsilon \in \R_{>0}$ there exists a time $T(\epsilon,\delta) \in \R_{>0}$ such that for all $\bar{v}_0 \notin \mc Z$
 \begin{align}\label{eq:Attractivity}
  \norm{ \bar{v}_0 }_{\mc T_{\bar{f},0}} < \delta \implies \norm{\varphi_{\bar{f}}(t,\bar{v}_0)}_{\mc T_{\bar{f},0}} < \epsilon \quad \forall t \geq T.
 \end{align}

According to Theorem \ref{th.convergence.to.S} the states $\bar{v}$ are exponentially stable with respect to $\mc S$ and we directly obtain the following property.\\
\emph{Property 1:} For every $\delta_{\mc S} \in \R_{>0}$ and every $\mu_{\mc S} \in \R_{>0}$ there exists a time $T_{\mu_{\mc S}}(\mu_{\mc S},\delta_{\mc S}) \in \R_{>0}$ such that for all $\bar{v}_0 \in \R^n \setminus \mc Z$:
 \begin{align}
  \norm{ \bar{v}_0 }_{\mc S} < \delta_{\mc S} \implies \norm{\varphi_{\bar{f}}(t,\bar{v}_0)}_{\mc S} < \mu_{\mc S} \quad \forall t \geq T_{\mu_{\mc S}}.
 \end{align}
Since $\mc T_{\bar{f},0} = \mc S \cap \mc A$ it remains to show that this property also holds for $\mc A$. To this end, we require the following properties which are a direct consequence of Proposition \ref{propo.wzeroinputconv}:

\noindent\emph{Property 2:}  For each compact subset ${\mathcal H} \subseteq \mc S\setminus \{0\}$ and each $\epsilon_{\mathcal H}$ there exists a time $T_{\epsilon_{\mathcal H}}({\mathcal H},\epsilon_{\mathcal H}) \in \R_{>0}$ such that for all $\bar{v}_0 \in {\mathcal H}$:
\begin{align}\label{eq.compact.attractivity}
    \norm{\varphi_{\bar{f}}(t,\bar{v}_0)}_{\mc A} < \epsilon_{\mathcal H},\quad  \forall t \geq T_{\epsilon_{\mathcal H}}.
\end{align}

\noindent\emph{Property 3:}
For each each $\epsilon \in \R_{>0}$ there exists a $\delta_{\mc A} \in \R_{>0}$ and a time $T_{\delta_{\mc A}} (\epsilon,\delta) \in \R_{>0}$ such that for all $\bar{v}_0 \in \mc S \setminus \{0\}$ it holds that:
\begin{align}\label{eq.prop.attract.stability}
  \norm{ \bar{v}_0 }_{\mc A} \leq \delta_{\mc A} \implies
  \begin{cases}
    \norm{\varphi_{\bar{f}}(t,\bar{v}_0)}_{\mc A} < \epsilon/2   & \forall t \in [0,T_{\delta_{\mc A}}],\\
    \norm{\varphi_{\bar{f}}(t,\bar{v}_0)}_{\mc A} < \delta_{\mc A}/2 & \forall t \geq T_{\delta_{\mc A}}.
  \end{cases}
\end{align}

We require the following preliminary result {that exploits the continuity of solutions of the vector field~\eqref{eq.desired.cl.magnitude.rotframe} with respect to the initial condition:}\\
\noindent\emph{Claim 1:} For each compact subset ${\mathcal H} \subseteq \mc S\setminus \{0\}$, each $T_b \in \R_{\geq 0}$, and each $\epsilon_b \in \R_{\geq 0}$, there is a $\mu({\mathcal H},T_b,\epsilon_b) \in \R_{>0}$ such that for each $\bar{v}^\prime_0 \in {\mathcal H}$ and each $\bar{v}_0\in \R^n\setminus\mc Z$ such that
\begin{align}\label{eq:distbound}
\begin{split}
\norm{\bar{v}_0 - \bar{v}^\prime_0} & \leq \mu  \implies \\
 &\norm{\varphi_{\bar{f}}(t,\bar{v}_0) - \varphi_{\bar{f}}(t,\bar{v}^\prime_0)} < \epsilon_b\,\, \forall t \in [0,T_b].
 \end{split}
\end{align}
\noindent\emph{Proof of Claim 1:}
The trajectories $\varphi_{\bar{f}}(t,\bar{v}_0)$ are well-defined for all times $t \in \R_{\geq 0}$ and all initial conditions $\bar{v}_0 \in \R^n$. Moreover, because $\bar{f}$ is continuously differentiable, the map $\phi_{\bar{f},t}: \bar{v}_0 \mapsto \varphi_{\bar{f}}(t,\bar{v}_0)$ is continuous for every $t \in \R_{>0}$. Suppose that ${\mathcal H}$, $T_b$, and $\epsilon_b$ are given, pick any $\bar{v}_0 \in {\mathcal H} \setminus \mc Z$. Because the map $\phi_{\bar{f},t}$ is continuous, there exists $\mu_{\bar{v}_0} \in \R_{>0}$ such that \eqref{eq:distbound} holds for all $\norm{\bar{v}_0 - \bar{v}^\prime_0}  \leq \mu_{\bar{v}_0}$. By compactness of ${\mathcal H}$ we define $\mu$ as the smallest such $\mu_{\bar{v}_0}$.\\

We now prove the Theorem. Pick any $\epsilon \in \mathbb{R}_{>0}$, any $\delta \in \mathbb{R}_{>0}$, and any $\bar{v}_0 \in \R^n \setminus \mc Z$ such that $\norm{\bar{v}_0}_{\mc T_{\bar{f},0}} \leq \delta$. By Proposition \ref{propo.zero.unstable} we have that $\varphi_{\bar{f}} (t,\bar{v}_0) \notin \mc Z$ for all $t \in \R_{\geq 0}$.  Next, consider the compact set
\begin{align}
\Delta(\delta) = \bigcup_{\norm{\bar{v}_0}_{\mc T_{\bar{f},0}} \leq \delta} \mc W_1(\bar{v}_0) \times \ldots \times \mc W_N(\bar{v}_0),
\end{align}
where the sets $\mc W_k$ are defined in~\eqref{eq.sets.wk}. By Proposition \ref{propo.wbound} the trajectories $\varphi_{\bar{f}} (t,\bar{v}_0)$ remain in $\Delta(\delta) \setminus \mc Z$ for all $t \in \R_{\geq 0}$, let ${\mathcal H}(\delta) \coloneqq \Delta(\delta) \cap (\mc S \setminus \{0\})$.\\

We can establish the following facts
\begin{enumerate}[(a)]
\item There exists a $\delta_{\mc A}<\epsilon/2$ and $T_{\delta_{\mc A}}(\delta_{\mc A},\epsilon) \in \R_{>0}$ that satisfy~\eqref{eq.prop.attract.stability}.
\item Let $\epsilon_{\mathcal H} = \delta_{\mc A}/2$, according to Property 2 we can pick $T_{\epsilon_{\mathcal H}}({\mathcal H},\delta_{\mc A}/2) \geq T_{\delta_{\mc A}}(\delta_{\mc A},\epsilon)$ such that \eqref{eq.compact.attractivity} is satisfied.
\item  Let $\epsilon_b = \delta_{\mc A}/2$ and $T_b = T_{\epsilon_{\mathcal H}}$, it follows from Claim 1 that there exists $\mu_{\mc S} \coloneqq \mu({\mathcal H},T_{\epsilon_{\mathcal H}}({\mathcal H},\delta_{\mc A}/2),\delta_{\mc A}/2)$ that satisfies~\eqref{eq:distbound}.
\item Given $\mu_{\mc S}$ as defined in (c), it follows from Property 1 that there exists $T_{\mu_{\mc S}}$ such that $\norm{\varphi_{\bar{f}} (t,\bar{v}_0)}_{\mc S} < \mu_{\mc S}$ for all $t>T_{\mu_{\mc S}}$. 
 \end{enumerate}

Since we know that after time $T_{\mu_{\mc S}}$ we are $\mu_{\mc S}$ close to $\mc S$, we can exploit the continuity arguments of Claim 1 and the stability properties of $\mc A$ in $\mc S$ to establish attractivity of $\mc T_{\bar{f},0} = \mc S\cap \mc A$. 

To this end, pick any $\bar{v}^\prime \in \mc S \setminus \{0\}$ such that $\norm{\varphi_{\bar{f}} (T_{\mu_{\mc S}},\bar{v}_0) - \bar{v}^\prime} \leq \mu_{\mc S}$. Such a $\bar{v}^\prime$ exists because $\norm{\varphi_{\bar{f}} (T_{\mu_{\mc S}},\bar{v}_0)} < \mu_{\mc S}$. By (b) we have that $\norm{\varphi_{\bar{f}}(T_{\epsilon_{\mathcal H}},\bar{v}^\prime)}_{\mc A} < \delta_{\mc A}/2$ and by (c) and (d) it follows that $\norm{\varphi_{\bar{f}}(T_{\epsilon_{\mathcal H}} +T_{\mu_{\mc S}} ,\bar{v}_0)}_{\mc A} < \delta_{\mc A}$. It remains to show that $\norm{\varphi_{\bar{f}}(t,\bar{v}_0)}_{\mc A} < \epsilon$ for all $t \geq T_{\epsilon_k} + T_{\mu_{\mc S}}$.

We have that $\norm{\varphi_{\bar{f}}(T_{\epsilon_{\mathcal H}} +T_{\mu_{\mc S}} ,\bar{v}_0)}_{\mc A} < \delta_{\mc A}$ and $\norm{\varphi_{\bar{f}}(T_{\epsilon_{\mathcal H}} +T_{\mu_{\mc S}} ,\bar{v}_0)}_{\mc S} < \mu_{\mc S}$. Using the same continuity argument we show that starting in a $\delta_{\mc A}$ neighborhood of $\mc A$ we cannot leave an $\epsilon$ neighborhood of $\mathcal A$.  Pick any $\bar{v}^\prime \in \mc S \setminus \{0\}$ such that $\norm{\varphi_{\bar{f}} (T_{\epsilon_{\mathcal H}} +T_{\mu_{\mc S}} ,\bar{v}_0) - \bar{v}^\prime} \leq \mu_{\mc S}$ and $\norm{\bar{v}^\prime}_{\mc A} < \delta_{\mc A}$. By (a) we have that $\norm{\varphi_{\bar{f}}(t,\bar{v}^\prime)}_{\mc A} < \epsilon/2$ for all $t \in [0, T_{\delta_{\mc A}}]$ and $\norm{\varphi_{\bar{f}}(T_{\delta_{\mc A}} ,\bar{v}^\prime)}_{\mc A} < \delta_{\mc A}/2$. By (c) and $T_{\epsilon_{\mathcal H}}({\mathcal H},\delta_{\mc A}/2) \geq T_{\delta_{\mc A}}(\delta_{\mc A},\epsilon)$ it follows that $\norm{\varphi_{\bar{f}}(t,\bar{v}_0)}_{\mc A} < \epsilon$ for all $t \in [T_{\epsilon_{\mathcal H}} +T_{\mu_{\mc S}}, T_{\epsilon_{\mathcal H}} +T_{\mu_{\mc S}}+T_{\delta_{\mc A}}]$ and $\norm{\varphi_{\bar{f}}(T_{\epsilon_{\mathcal H}} +T_{\mu_{\mc S}}+T_{\delta_{\mc A}},\bar{v}_0)}_{\mc A} < \delta_{\mc A}$. By induction it follows that $\norm{\varphi_{\bar{f}}(t,\bar{v}_0)}_{\mc A} < \epsilon$ for all $t > T_{\epsilon_{\mathcal H}} +T_{\mu_{\mc S}}$.

It follows that for all initial conditions $\bar{v}_0 \in \R^n \setminus \mc Z$ and every $\epsilon \in \R_{>0}$, \eqref{eq:Attractivity} holds for $T = T_{\epsilon_{\mathcal H}} +T_{\mu_{\mc S}}$. By Proposition \ref{propo.zero.unstable} the region of attraction $\mc Z$ of the equilibrium $\bar{v}=0$ is a zero measure set, i.e. the convergence property holds for almost every initial condition.
\end{IEEEproof}

\subsubsection{Proof of stability}
In Theorem~\ref{thm.sync}, we established that $\mc T_{\bar{f},0}$ is (almost) globally attractive. In order to show (almost) global asymptotic stability we need to show that the set $\mc T_{\bar{f},0}$ is also stable in the sense of Definition~\ref{def:stability}. To establish stability of the closed-loop dynamics, we define the function $\gamma_{\mc S}(\gamma_{\mc A})=\frac{v^\star_{\min}}{v^\star_{\max}} \frac{\gamma_{\mc A} (v^\star_{\min}   - \gamma_{\mc A})}{\bar{\sigma}(\mc K - \mc L)}$ and the set
  \begin{align*}
   \mc M(\gamma_{\mc A}) \coloneqq \!\left\{ \bar{v} \in \mathbb{R}^n \left\vert\!\!\!
   \begin{array}{l}
    \norm{\bar{v}}_{\mc S} \leq \gamma_{\mc S}(\gamma_{\mc A}), \\
    \norm{\bar{v}_k}_{\mc A_k} \leq  \gamma_{\mc A} \qquad \forall k\in\mc V
   \end{array}\!\!\! \right. \right\}\!.
\end{align*}
\begin{proposition}{\bf(Invariant neighborhood of $\boldsymbol{\mc T_{\bar{f},0}}$)}\label{propo.invset}\\
Under Condition~\ref{cond.small.angle}, the set $\mc M(\gamma_{\mc A})$ is non-empty and invariant with respect to the dynamics \eqref{eq.desired.cl.magnitude.rotframe} for all $\gamma_{\mc A} \in \mathbb{R}_{[0,\frac{1}{2}{v^\star_{\min}})}$. For any $\gamma_{\mc A,1},\gamma_{\mc A,2} \in \mathbb{R}_{[0,\frac{1}{2}v^\star)}$ such that $\gamma_{\mc A,2} < \gamma_{\mc A,1}$ it holds that $\mc M(\gamma_{\mc A,2}) \subset \mc M(\gamma_{\mc A,1})$. Moreover, for any $\gamma_{\mc A} \in \mathbb{R}_{(0,\frac{1}{2} {v^\star_{\min}})}$ the set $\mc T_{\bar{f},0}$ lies in the interior of $\mc M(\gamma_{\mc A})$ and the interior of $\mc M(\gamma_{\mc A})$ is non-empty.
\end{proposition}
The proof is provided in the Appendix. Broadly speaking, the proof establishes that trajectories stay close to $\mc A$ if they start sufficiently close to $\mc S$ and $\mc A$. In combination with exponential stability of $\mc S$ this allows to construct a family of invariant sets around $T_{f,\omega_0}$. The next Theorem uses these invariant sets to establish stability of the target set $\mc T_{f,\omega_0}$.

\begin{theorem}{\bf (Stability of $\boldsymbol{\mc T_{f,\omega_0}}$)}\label{thm.stab}\\
 Under Condition~\ref{cond.small.angle}, the dynamical system~\eqref{eq.desired.cl.magnitude} is stable with respect to the set $\mc T_{f,\omega_0}$.
 \end{theorem}
\begin{IEEEproof}
By Proposition~\ref{propo.equivalent.rotating}, the Theorem's statement is equivalent to stability of \eqref{eq.desired.cl.magnitude.rotframe} with respect to the set $\mc T_{\bar{f},0}$ in the rotating coordinate frame. According to Proposition \ref{propo.invset} for all $\gamma_{\mc A} \in \mathbb{R}_{(0,\frac{1}{2} {v^\star_{\min}})}$ the set $\mc M(\gamma_{\mc A})$ is invariant with respect to the dynamics \eqref{eq.desired.cl.magnitude.rotframe} and contains $\mc T_{\bar{f},0}$ in its interior. Because $\mc M(\gamma_{\mc A})$ is compact for all $\gamma_{\mc A} \in \mathbb{R}_{(0,\frac{1}{2} {v^\star_{\min}})}$ it follows that for every $\epsilon \in \mathbb{R}_{>0}$ there exists a $\gamma_\epsilon \in \mathbb{R}_{(0,\frac{1}{2}v^\star_{\min})}$ such that $\mc M(\gamma_\epsilon) \subset \mathscr{B}_{\mc T_{\bar{f},0}}(\epsilon)$, i.e., $\mc M(\gamma_\epsilon)$ is contained in an $\epsilon$-neighborhood of $\mc T_{\bar{f},0}$. Next, because $\mc M(\gamma_{\mc A})$ has non-empty interior for all $\gamma_{\mc A} \in \mathbb{R}_{(0,\frac{1}{2} {v^\star_{\min}})}$ it follows that there exists a corresponding $\delta \in \mathbb{R}_{>0}$ such that $\mathscr{B}_{\mc T_{\bar{f},0}}(\delta) \subset \mc M(\gamma_\epsilon) \subset \mathscr{B}_{\mc T_{\bar{f},0}}(\epsilon)$. Because $\mc M(\gamma_\epsilon)$ is invariant, it follows that $\bar{v}_0 \in \mathscr{B}_{\mc T_{\bar{f},0}}(\delta)$ implies that $\varphi_{\bar{f}}(t,\bar{v}_0) \in \mc M(\gamma_\epsilon) \subset \mathscr{B}_{\mc T_{\bar{f},0}}(\epsilon)$ holds for all $t \geq 0$. This is precisely the definition of stability with respect to the set $\mc T_{\bar{f},0}$.
\end{IEEEproof}

Finally, using the results developed in this section, we can provide the proof of Theorem~\ref{thm.ags}.

\vspace{1em}

\emph{Proof of Theorem~\ref{thm.ags}{:}}
The Proof directly follows from Definition \ref{def:ags}, Theorem \ref{thm.sync}, and Theorem \ref{thm.stab}.  \hfill \IEEEQED

\section{Power systems case study}\label{sec.example}
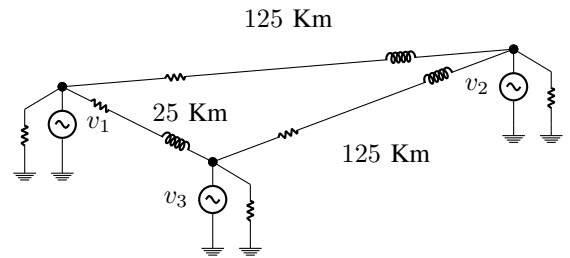
\begin{figure}[htbp]
\begin{center}
\begin{circuitikz}[american voltages]

\ctikzset{bipoles/resistor/height=0.15}
\ctikzset{bipoles/resistor/width=0.4}

\ctikzset{bipoles/generic/height=0.15}
\ctikzset{bipoles/generic/width=0.4}

\ctikzset{bipoles/length=.6cm}

\coordinate (I1) at (0,0);
\coordinate (I2) at (6,0.5);
\coordinate (I3) at (2,-1);

\draw ($ (I1) + (0,-1) $) node [ground] (g1) {};
\draw (g1) to[sV, l_=$v_1$] (I1);
\draw ($ (I1) + (-0.5,-1) $) node [ground] (g12) {};
\draw (g12) to[R] ($ (I1) + (-0.5,-0.3) $);
\draw (I1) to[short] ($ (I1) + (-0.5,-0.3) $);

\draw ($ (I2) + (0,-1) $) node [ground] (g2) {};
\draw (g2) to[sV, l=$v_2$] (I2);
\draw ($ (I2) + (0.5,-1) $) node [ground] (g22) {};
\draw (g22) to[R] ($ (I2) + (0.5,-0.3) $);
\draw (I2) to[short] ($ (I2) + (0.5,-0.3) $);

\draw ($ (I3) + (0,-1) $) node [ground] (g3) {};
\draw (g3) to[sV, l=$v_3$] (I3);
\draw ($ (I3) + (0.5,-1) $) node [ground] (g32) {};
\draw (g32) to[R] ($ (I3) + (0.5,-0.3) $);
\draw (I3) to[short] ($ (I3) + (0.5,-0.3) $);

\node (I12) at ($(I1)!0.5!(I2)$) [label={[xshift=0cm, yshift=0.3cm]$125~\text{Km}$}]{};
\draw 		(I1) 
to[R,*-] (I12)
to [L,-*] (I2);

\node (I13) at ($(I1)!0.5!(I3)$)[label={[xshift=0.7cm, yshift=-0.2cm]$25~\text{Km}$}]{};
\draw 		(I1) 
to[R,*-] (I13)
to [L,-*] (I3);

\node (I32) at ($(I3)!0.5!(I2)$) [label={[xshift=0.3cm, yshift=-1.0cm]$125~\text{Km}$}]{};
\draw 		(I3) 
to[R,*-] (I32)
to [L,-*] (I2);

 \end{circuitikz}
\caption{Inverter based grid}
\label{fig.inverter.example}
\end{center}
\end{figure}

In this section, we {illustrate the theoretical results with a simplified power system case study}. We consider three inverters connected by the grid as illustrated in Figure~\ref{fig.inverter.example}. The grid base power is 1 GW, the base voltage {is} 320 {kV,} and the base frequency {is} 50 Hz. The line resistance is 0.03 Ohm/Km, and the line reactance is 0.3 Ohm/km (at the nominal frequency). Therefore, the reactance/resistance ratio $\omega_0\rho$ of the transmission lines is given by $\omega_0\rho=10$, and the angle $\kappa$ can be computed as $\kappa = \tan^{-1}(10)=84.2894^\circ$.

\setlength{\tabcolsep}{4pt}
\renewcommand{\arraystretch}{1.2}
\begin{table}
\caption{Power and voltage set-points for the three inverters}
\label{tab.pf}
\begin{center}
\begin{tabular}{ @{} c|rrr|rrr @{} } 
\bottomrule
\multicolumn{6}{>{\columncolor[gray]{.95}}l}{Power and voltage set-points during the simulation.} \\
\toprule
 &  \multicolumn{2}{r}{$t=5\, s~~$} & & \multicolumn{2}{r}{$t=10\,s~~$} \\

 $k$  & $p_k^\star$ p.u. & $q_k^\star$ p.u.  & $v_k^\star$ p.u.  & $p_k^\star$ p.u. & $q_k^\star$ p.u.  & $v_k^\star$ p.u.  \\[0.1em]
\toprule 
 1 & $0.1458$ & $0.0432$ & ${1.01}$  & $0.1458$ & $0.0432$ & ${1.01}$ \\
 2 & $0.7066$  & $-0.0793$ & $1$ & $0.7066$  & $-0.0793$ & $1$ \\
 3  & $-0.8509$ &  $0.0803$ & $1$ & $\mathbf{-0.3509}$ &  $0.0803$ & $1$ \\
\toprule
\bottomrule

\end{tabular}
\end{center}
\end{table}
\begin{figure}[htbp]
\begin{center}
\input{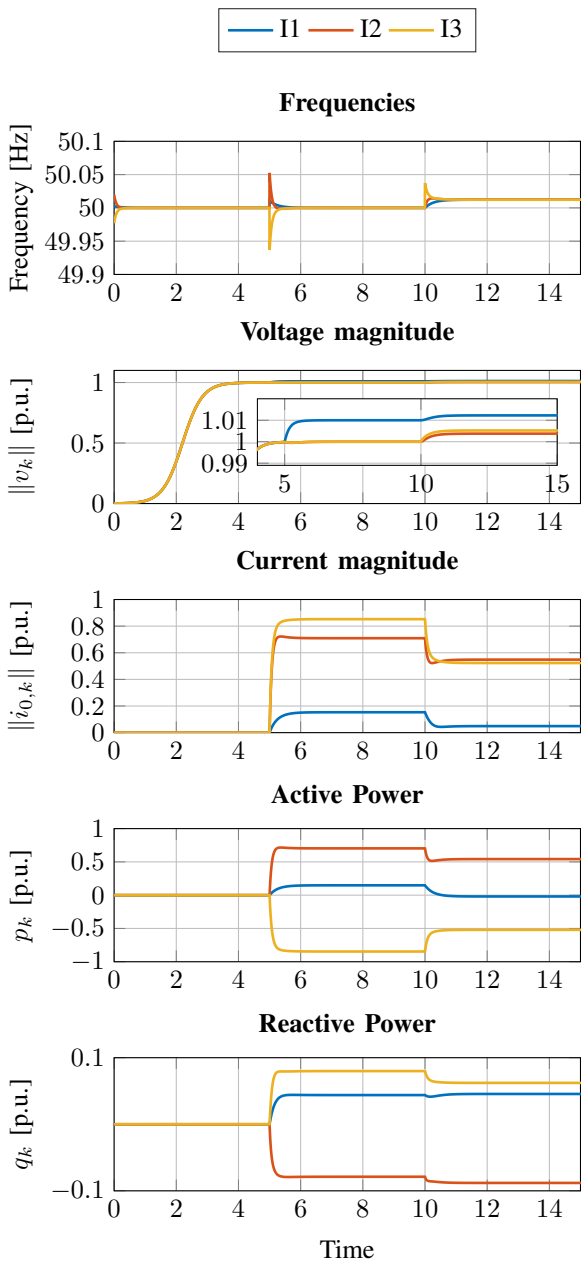}
\caption{ At $t=0\,s$, the grid is black started giving the inverters zero power set-points. At $t=5\,s$ the grid converges to the set-points provided in the first column of Table~\ref{tab.pf}. At $t=10\,s$, inverter 3 unilaterally changes its active power set-point $p_1^\star$ to $-0.3509$. To compensate, {Inverters} 1 and 2 reduce their power injection and the systems settles to a different synchronous solution.}
\label{fig.blackstart}
\end{center}
\end{figure}
We compute controllers of the form~\eqref{eq.controller} using~\eqref{eq.K.in.PQ}, i.e.,  the matrices $K_k$ for each inverter are constructed only using the power and voltage magnitude set-points. We simulate a black-start and two dispatches: the first dispatch assigns a set of feasible set-points to each inverter, the second dispatch assigns a set of set-points corresponding to an infeasible power-flow. In Figure~\ref{fig.blackstart} we show the frequency, the voltage magnitudes, current magnitudes, and the active and reactive powers injections of each inverter.
At $t=0\,s$ we simulate a black-start for the inverter based grid. We initialize the voltages to $v_k(0) = 10^{-3} [1,\, 1]^\top$, for $k=1,2,3$ and we provide the inverters with zero power set-points and voltage set-points $v_k^\star=1$ p.u. At $t=5\,s$ we simulate a dispatch to the set-points presented in the respective column of Table~\ref{tab.pf}. The power-flow solution in Table~\ref{tab.pf} provides the following angles at steady state $\theta_{21}^\star = 0^\circ$ and  $\theta_{31}^\star = -3^\circ$.  We chose $\eta=0.0015$ and $\alpha=0.01$ so that condition~\eqref{eq.condtion.js} is satisfied. As expected from Theorem~\ref{thm.ags}, the inverters synchronize to the correct {relative phase} and frequency, {and} reach the prescribed voltage magnitudes and power injections. Then, at $t=10\,s$, inverter 3, which in this simulation is acting as a load (i.e., has a negative active power set-point), unilaterally steps up it's active power set-point $p_1^\star$ by $0.5$ p.u., while all other inverters keep their set-points and controllers unchanged. Therefore, the set-points no longer correspond to a feasible power-flow and guarantees given by Theorem~\ref{thm.ags} are lost. We note, however, that the system remains stable and synchronous and, as expected from the droop-like behavior (see Proposition \ref{propo.droop.like}) of the closed loop, inverter 2 and 3 share the power step down needed to accommodate for the reduced load.

\section{Conclusion and outlook}\label{sec.conclusion}
In this paper, we propose{d} a new decentralized control method for coupled oscillators with (almost) global {asymptotic stability} guarantees and we explore{d} its application to the control of grid-connected power inverters. {The proposed methods ensure} synchronization and almost global {stability} of a solution with the desired power injections. {Some of the technical assumptions needed to prove the results are restrictive, in particular in the context of applying the proposed control scheme to the power system.} This opens numerous further research directions. {Among other assumptions,} in this paper, we consider{ed} a network containing only inverters acting as source or as a load, and we neglect{ed} the dynamics of the transmission lines. A question that arises naturally is the analysis of the proposed control method for the case of {heterogeneous} dynamic lines, structure preserving models with load buses, and the presence of synchronous generators. {Finally, the robustness properties of the proposed control scheme with respect to inconsistent power or magnitude set-points need to \red{be} analyzed.} From the droop-like behavior and preliminary simulation results such as the one presented in Section~\ref{sec.example}, we expect the system to be well-behaved and exhibit beneficial power-sharing characteristics. \red{However, thus far we lack analytical conditions that guarantee both stability and robustness. Finally, experimental validation of the proposed controller is the subject of ongoing work and preliminary results can be found in \cite{SC+19}.}

\appendices
\section{Proofs of technical results}
The appendix contains the proofs of some technical results from Sections~\ref{sec.consensus.control} and~\ref{sec.proofs}.

\emph{Proof of Proposition~\ref{propo.implementable}:}~
Since $\kappa = \tan^{-1}(\rho\,\omega_0)$, and using the trigonometric identities
\begin{align}\label{eq.trig}
\sin(\kappa) &=~ \frac{\rho\,\omega_0}{\sqrt{1+(\rho\,\omega_0)^2}} = \frac{\ell_{jk}\,\omega_0}{\sqrt{r_{jk}^2+\omega_0^2\ell_{jk}^2}}\\
\cos(\kappa) &=~ \frac{1}{\sqrt{1+(\rho\,\omega_0)^2}} = \frac{r_{jk}}{\sqrt{r_{jk}^2+\omega_0^2\ell_{jk}^2}},
\end{align}
we can write~\eqref{eq.powers} as
\begin{align}
\begin{split}\label{eq.powers.ref}
p^\star_k ={v_k^{\star2}} \sum_{j:(j,k)\in\mc E} w_{jk}  \left[ \cos(\kappa)  - \frac{v_j^\star }{v_k^\star }\cos(\theta_{jk}^\star -\kappa) \right]\\
q^\star_k ={v_k^{\star2} }\sum_{j:(j,k)\in\mc E} w_{jk} \left[  \sin(\kappa) +  \frac{v_j^\star }{v_k^\star } \sin(\theta_{jk}^\star -\kappa) \right].
\end{split}
\end{align}
{From}~\eqref{eq.powers.ref} we note that
{
\begin{align*}
\begin{bmatrix}
p_k^\star & q_k^\star\\
-q_k^\star & p_k^\star
\end{bmatrix} = 
v_k^{\star2} R(\kappa)^\top  \sum\nolimits_{j:(j,k)\in\mc E}  w_{jk}(I_2-\frac{v^\star_j}{v^\star_k} R(\theta_{jk}^\star)).
\end{align*}
and \eqref{eq.K.in.PQ} immediately follows from \eqref{eq.Kk.def}.} \hfill \IEEEQED


\emph{Proof of Proposition~\ref{propo.droop.like}}~
Let us express each vector $v_k=[v_\alpha,v_\beta]^\top$ as a complex number $\mathbf v_k = v_{i,\alpha} + j v_{i,\beta}$. 
We define by $y_{ik}\coloneqq (\omega_0\ell_{ik})^{-1}$, then, considering the expression for $\mc K$ given in~\eqref{eq.K.in.PQ} we can re-write~\eqref{eq.desired.cl.magnitude} as
\begin{align}
\begin{split}\label{eq.complex.valued.dyn}
&\dot \nu_k \mathrm e ^{j\theta_k} + j \dot \theta_k\nu_k \mathrm e ^{j\theta_k} = j \left(\omega_0+\eta \frac{p^\star_k}{{v^{\star2}_k}}\right)\nu_k \mathrm e ^{j\theta_k}  \\
&+ \eta \sum_{i:(i,k)\in\mc E}\,y_{ik} \nu_i \mathrm e ^{j(\theta_i)} -  \eta \sum_{i:(i,k)\in\mc E}\,y_{ik} \nu_k \mathrm e ^{j(\theta_k)}\\
& + \eta\frac{q^\star_k}{{v^{\star2}_k}} \nu_k \mathrm e ^{j\theta_k} + \frac{\alpha}{{v^{\star}_k}}({v^{\star}_k}-\nu_k)\nu_ke ^{j\theta_k}.
\end{split}
\end{align} 
After separating the real and imaginary part of~\eqref{eq.complex.valued.dyn} and performing simple trigonometric manipulations we obtain:
\begin{align*}
\dot \nu_k = &~ \eta\left( \sum_{i:(i,k)\in\mc E} \,y_{ik}\left( \frac{\nu_i\nu_k}{\nu_k^2} \cos(\theta_i-\theta_k) -1\right){+}\frac{q^\star_k}{{v^{\star2}_k}}\right) \nu_k\\
&+~\frac{\alpha}{{v^\star_k}}({v^{\star}_k}-\nu_k)\nu_k\\
\dot \theta_k = &~  \omega_0+\eta \left( \frac{p_k^\star}{{v^{\star2}_k}} -  \sum_{i:(i,k)\in\mc E} \,y_{ik} \frac{\nu_i\nu_k}{\nu_k^2} \sin(\theta_i-\theta_k) \right).
\end{align*}
{Using \eqref{eq.trig} and similar reformulations as in the proof of Proposition~\ref{propo.implementable} we obtain \eqref{eq.droop.like}.} \hfill\IEEEQED

\emph{Proof of Proposition~\ref{propo.zero.unstable}}~
  The proof is based on linearizing the dynamics \eqref{eq.desired.cl.magnitude.rotframe} at $\bar{v}=0$. The Jacobian of $\Phi_k(\bar{v}_k) \bar{v}_k$ with respect to $\bar{v}_k$ is given by:
  \begin{align*}
   \frac{\partial}{\partial \bar{v}_k}  \Phi_k(\bar{v}_k) \bar{v}_k = \begin{cases} I_2 - \frac{1}{{v^\star_k} } \left(\frac{\bar{v}_k \bar{v}^{\top}_k}{\norm{\bar{v}_k}} + I_2 \norm{\bar{v}_k }\right) & \bar{v}_k \neq 0 \\
                                                      I_2 & \bar{v}_k = 0
                                    \end{cases}.
  \end{align*}
  Moreover, considering the limit $\lim_{\bar{v}_k \to 0} \frac{\bar{v}_k \bar{v}^{\top}_k}{{v^\star_k} \norm{\bar{v}_k}} = 0$ it can be verified that the derivative is continuous {at $\bar{v}$} and $\Phi(\bar{v})\bar{v}$ is continuously differentiable. The Jacobian of $\bar{f}(\bar{v})$ at $\bar{v}=0$ is given by $\frac{\partial \bar{f}(\bar{v})}{\partial \bar{v}}\vert_{\bar{v}=0} =\eta (\mc K -\mc L) + \alpha I_{2N}$.  Theorem {\ref{th.convergence.to.S} implies that the set $\mc S$ is positively invariant under the dynamics \eqref{eq.desired.cl.magnitude.rotframe}. Therefore} $\mc S \subseteq \ker(\mc K -\mc L)$, and it follows that at least two eigenvalues of $(\mc K -\mc L)$ are zero. Because of $\alpha >0$, the real part of at least two eigenvalues of the Jacobian $\frac{\partial \bar{f}(\bar{v})}{\partial \bar{v}}\vert_{\bar{v}=0}$ is positive. Therefore, $\bar{v}=0$ is an unstable equilibrium. By Proposition 11 of \cite{Monzon2006} this, and the fact that $\bar{f}(\bar{v})$ is continuously differentiable, implies that the set {$\mc Z$} of initial conditions $\bar{v}_0$ such that $\varphi_{\bar{f}}(t,\bar{v}_0) \to 0$ as $t \to \infty$ has zero Lebesgue measure.  \hfill\IEEEQED

\emph{Proof of Proposition~\ref{propo.wbound}}~
Consider the Lyapunov-like functions $W_k(\bar{v}_k) = \Phi_k(\bar{v}_k)^2$ for every $k \in \mc V$. The time derivative of $W_k$ is given by
  \begin{align}\begin{split}\label{eq.wk.decr}
\frac{\diff}{\diff t} W_k &= -2 W_k \norm{\bar{v}_k} - \frac{2}{{v^\star_k}} \Phi_k(\bar{v}_k) \frac{\bar{v}^\top_k}{\norm{\bar{v}_k}} e_{\theta,k} \\
&\leq -2 W_k \norm{\bar{v}_k} + \frac{2}{{v^\star_k}} \sqrt{W_k} \bar{\sigma}(\mc K - \mc L) \norm{\bar{v}}_{\mc S},
  \end{split}\end{align}
{where $e_{\theta,k} {\in \mathbb{R}^2}$ is the $k^{\text{th}}$ component of ${e_{\theta}=(\mc K - \mc L)\bar v \in \mathbb{R}^{2N}}$} and  $\bar{\sigma}(\mc K - \mc L)$ denotes the largest singular value of the matrix $\mc K - \mc L$.

By Theorem \ref{th.convergence.to.S} it holds that $\norm{\varphi_{\bar{f}}(t,\bar{v}_0)}_{\mc S} \leq \norm{\bar{v}_0}_{\mc S}$ for all $t \geq 0$. Moreover, it can be verified that $W_k>1$ implies $\norm{\bar{v}_k} > 2 {v^\star_k}$, i.e., for all $W_k>1$ and all $t \geq 0$ it holds that:
\begin{align*}
\frac{\diff}{\diff t} W_k &\leq -4 {v^\star_k} W_k + \frac{2}{{v^\star_k}} \sqrt{W_k} \bar{\sigma}(\mc K - \mc L) \norm{\bar{v}(0)}_{\mc S}\\
&{=} -4 \sqrt{W_k} \left({v^\star_k} \sqrt{W_k} - \frac{1}{2 {v^\star_k}} \bar{\sigma}(\mc K - \mc L) \norm{\bar{v}(0)}_{\mc S}\right).
\end{align*}
It directly follows that $\frac{\diff}{\diff t} W_k < 0$ if $\sqrt{W_k} > 1$ and
$\sqrt{W_k} > \frac{1}{2 {v^{\star2}_k}} \bar{\sigma}(\mc K - \mc L) \norm{\bar{v}(0)}_{\mc S}$. This implies
that for every $\bar{v}_{k,0}$ there exists a constant $\hat{W}_k \geq W_k(\bar{v}_{k,0})$ such that 
$W_k \leq \hat{W}_k$ holds for all $t \geq 0$. \hfill\IEEEQED

\emph{Proof of Proposition~\ref{propo.invset}:}~{
Let us recall the function {$\gamma_{\mc S}(\gamma_{\mc A})=\frac{v^\star_{\min}}{v^\star_{\max}} \frac{\gamma_{\mc A} (v^\star_{\min}   - \gamma_{\mc A})}{\bar{\sigma}(\mc K - \mc L)}$}.
By definition it holds that $\gamma_{\mc S}(0)=0$, and it can be verified that $\gamma_{\mc S}(\gamma_{\mc A}) \in \mathbb{R}_{>0}$ and $\tfrac{\diff \gamma_{\mc S}}{\diff \gamma_{\mc A}} \in \mathbb{R}_{>0}$ {holds for all} $\gamma_{\mc A} \in \mathbb{R}_{[0,\frac{1}{2}{v^\star_{\min}})}$, i.e. the function $\gamma_{\mc S}(\gamma_{\mc A})$ is strictly positive and strictly increasing for $\gamma_{\mc A} \in \mathbb{R}_{[0,\frac{1}{2}{v^\star_{\min}})}$. It directly follows from Theorem \ref{th.convergence.to.S} that, for every fixed $\gamma_{\mc A} \in \mathbb{R}_{[0,\frac{1}{2}{v^\star_{\min}})}$, $\norm{\bar{v}_0}_{\mc S} \leq \gamma_{\mc S}(\gamma_{\mc A})\implies\norm{\varphi_{\bar{f}}(t,\bar{v}_0)}_{\mc S} \leq \gamma_{\mc S}(\gamma_{\mc A})$ for all $t \in \mathbb{R}_{\geq 0}$, i.e., the set $\{ \bar{v} \in \mathbb{R}^n \; \vert\; \norm{\bar{v}}_{\mc S} \leq \gamma_{\mc S}(\gamma_{\mc A}) \}$ is non-empty and invariant for all $\gamma_{\mc A} \in \mathbb{R}_{[0,\frac{1}{2}{v^\star_{\min}})}$.

Considering \eqref{eq.wk.decr}  and letting $\zeta \coloneqq \frac{1}{{v^\star_{\min}}}  \bar{\sigma}(\mc K - \mc L)$ it holds for all $\norm{\bar{v}}_{\mc S} \leq \gamma_{\mc S}(\gamma_{\mc A})$ and all $\gamma_{\mc A} \in \mathbb{R}_{[0,\frac{1}{2}{v^\star_{\min}})}$ that
\begin{align}
\ddt W_k \leq -2 \sqrt{W_k} \left( \sqrt{W_k} \norm{\bar{v}_k} - \zeta \gamma_{\mc S}(\gamma_{\mc A}) \right).
\end{align}
Next, it follows from $\norm{\bar{v}_k}^2_{\mc A_k} = \big(\norm{\bar{v}_k} - {v^\star_k}\big)^2$ that $\norm{\bar{v}_k} \geq {v^\star_k}  - \gamma_{\mc A}$ holds for all $\bar{v}_k$ such that $\norm{\bar{v}_k}_{\mc A_k} \leq \gamma_{\mc A}$. Moreover, we can use $v^\star_k \sqrt{W_k} = \norm{\bar{v}_k}_{\mc A_k}$ to obtain the following bound for all $\bar{v}$ on the boundary of $\mc M(\gamma_{\mc A})$, i.e., for $\norm{\bar{v}_k}_{\mc A_k} = \gamma_{\mc A}$, all $\gamma_{\mc A} \in \mathbb{R}_{[0,\frac{1}{2}{v^\star_{\min}})}$, and all $k\in\mc V$:
\begin{align}
 \ddt W_k \leq -2 \sqrt{W_k} \left( \left(1  - \frac{\gamma_{\mc A}}{{v^\star_k}}\right) \gamma_{\mc A} - \zeta \gamma_{\mc S}(\gamma_{\mc A}) \right).
\end{align}
It can be verified that $\ddt W_k \leq 0$ holds for all $\bar{v}$ on the boundary of $\mc M(\gamma_{\mc A})$ if the following condition holds
\begin{align}\label{eq:Wbound}
 \gamma_{\mc A} ({v^\star_{\min}}   - \gamma_{\mc A})  \geq {v^\star_{\max}} \zeta \gamma_{\mc S}(\gamma_{\mc A}).
\end{align}
This condition is satisfied by definition of $\gamma_{\mc S}(\gamma_{\mc A})$ and implies that the set $\{ \bar{v} \in \mathbb{R}^n \; \vert\; \norm{\bar{v}_k}_{\mc A_k} \leq  \gamma_{\mc A} \}$ is non-empty and invariant with respect to the dynamics \eqref{eq.desired.cl.magnitude.rotframe} for all $\gamma_{\mc A} \in \mathbb{R}_{[0,\frac{1}{2}{v^\star_{\min}})}$. Because the intersection of invariant sets is invariant \cite{BM08} and $\{ \bar{v} \in \mathbb{R}^n \; \vert\; \norm{\bar{v}}_{\mc S} \leq \gamma_{\mc S}(\gamma_{\mc A}) \}$ is invariant for $\gamma_{\mc A} \in \mathbb{R}_{[0,\frac{1}{2}{v^\star_{\min}})}$ it follows that $\mc M(\gamma_{\mc A})$ is invariant for all $\gamma_{\mc A} \in \mathbb{R}_{[0,\frac{1}{2}{v^\star_{\min}})}$. Moreover, because of $\gamma_{\mc S}(0)=0$ it holds that $\mc M (0) = \mc T_{\bar{f},0}$, i.e., $\mc M (0)$ is non-empty and it immediately follows that $\mc M(\gamma_{\mc A})$ contains a neighborhood of $\mc T_{\bar{f},0}$ for all $\gamma_{\mc A} \in \mathbb{R}_{(0,\frac{1}{2} {v^\star_{\min}})}$. In other words, for any $\gamma_{\mc A} \in \mathbb{R}_{(0,\frac{1}{2} {v^\star_{\min}})}$ the set $\mc T_{\bar{f},0}$ lies in the interior of ${\mc M(\gamma_{\mc A})}$ and the interior of ${\mc M(\gamma_{\mc A})}$ is non-empty.

Finally, consider any $\gamma_{\mc A,1},\gamma_{\mc A,2} \in \mathbb{R}_{[0,\frac{1}{2}v^\star)}$ such that $\gamma_{\mc A,2} < \gamma_{\mc A,1}$. Because $\gamma_{\mc S}(\gamma_{\mc A})$ is strictly increasing for $\gamma_{\mc A} \in \mathbb{R}_{[0,\frac{1}{2}{v^\star_{\min}})}$ it holds that $\gamma_{\mc S}(\gamma_{\mc A,2}) < \gamma_{\mc S}(\gamma_{\mc A,1})$ and, by construction of $\mc M(\gamma_{\mc A})$ it directly follows that $\mc M(\gamma_{\mc A,2}) \subset \mc M(\gamma_{\mc A,1})$.}  \hfill\IEEEQED

\section{Clarke Transformation}\label{app.clarke}
{
Given a three phase signal $v_{abc} = (v_a,v_b,v_c)${,} we define the amplitude preserving Clarke transformation~\cite{clarke1943circuit} as 
\[
v_{{\alpha \beta \gamma }}(t)=\begin{bmatrix}v_{\alpha}(t)\\v_{\beta}(t)\\ v_{\gamma}(t) \end{bmatrix}={\frac  {2}{3}}{\begin{bmatrix}1&-{\frac  {1}{2}}&-{\frac  {1}{2}}\\0&{\frac  {{\sqrt  {3}}}{2}}&-{\frac  {{\sqrt  {3}}}{2}}\\{\frac  {1}{2}}&{\frac  {1}{2}}&{\frac  {1}{2}}\\\end{bmatrix}}{\begin{bmatrix}v_{a}(t)\\v_{b}(t)\\v_{c}(t)\end{bmatrix}}.
\]
Note that if $v_{abc}$ is \emph{balanced} {(i.e., $v_a+v_b+v_c = 0$)}, then  
\[
v_{{\alpha \beta \gamma }}(t)=\begin{bmatrix}v_{\alpha}(t)\\v_{\beta}(t)\\0 \end{bmatrix}={\frac  {2}{3}}{\begin{bmatrix}1&-{\frac  {1}{2}}&-{\frac  {1}{2}}\\0&{\frac  {{\sqrt  {3}}}{2}}&-{\frac  {{\sqrt  {3}}}{2}}\\{\frac  {1}{2}}&{\frac  {1}{2}}&{\frac  {1}{2}}\\\end{bmatrix}}{\begin{bmatrix}v_{a}(t)\\v_{b}(t)\\v_{c}(t)\end{bmatrix}}.
\]
Hence, the \emph{two-dimensional} signal $v_{\alpha \beta}:= (v_\alpha,v_\beta)$ contains the same information as the original three-phase signal $v_{abc}$.
}

\bibliographystyle{IEEEtran}
\bibliography{IEEEabrv,bib_file}

\end{document}